\documentstyle[11pt]{article}
\input amssym.def
\input amssym

\newtheorem{theorem}{Theorem}[section]
\newtheorem{lemma}[theorem]{Lemma}
\newtheorem{proposition}[theorem]{Proposition}

\newtheorem{conjecture}[theorem]{Conjecture}

\begin{document}

\begin{titlepage}
\title{\bf Volumes of hyperbolic manifolds and mixed Tate motives}

\author{Alexander Goncharov}
\end{titlepage}
\maketitle
\tableofcontents

\def\QQ{{\Bbb Q}}
\def\CC{{\Bbb C}}
\def\RR{{\Bbb R}}%
\begin{displaymath}
\end{displaymath}

\def\ZZ{{\Bbb Z}}
\def\bqq{{\bar \QQ}}
\def\cl{{\Cal L}}
\def\car{{\Cal R}}
\def\cb{{\Cal B}}
\def\cd{{\cal D}}
\def\ch{{\Cal H}}

\def\g{\gamma}
\def\De{\Delta}
\def\a{\alpha}
\def\vp{{\varphi}}
\def\b{\beta}
\def\G{\Gamma}
\def\de{\delta}
\def\La{\Lambda}
\def\ve{\varepsilon}

\def\vo{{\hbox{vol}}}
\def\re{{\hbox{Re}}}
\def\im{{\hbox{Im}}}
\def\ot{{\otimes}}
\def\cone{{\hbox{Cone}}}
\def\tro{{\tilde \rho}}
\def\tq{{\tilde Q}}
\def\tc{{\tilde c}}

\def\ts{{\tilde S}}
\def\tv{{\tilde v}}
\def\hx{{\hat x}}
\def\hz{{\hat z}}
\def\hc{{\hat c}}
\def\hg{{\hat g}}
\def\wh{\widehat}



\section { Introduction}

By hyperbolic manifold we will mean an orientable complete Riemannian manifold
with constant sectional curvature -1.
\vskip 3mm \noindent
 {\bf 1. Volumes of hyperbolic $2m$-manifolds}. Let $M^n$ be an
$n$-dimensional hyperbolic manifold with finite volume $\vo({M^n})$.  Suppose
that $n=2m$ is an even number.  Then according to the Gauss-Bonnet theorem
(see, for example, [Ch])
$$
\vo(M^{2m}) = -c_{2m} \cdot \chi(M^{2m})
$$
where $c_{2m} = 1/2 \times $(volume of sphere $S^{2m}$ of radius 1) and $\chi
(M^{2m})$ is the Euler characteristic of $M^{2m}$.  This is straightforward
for compact manifolds and a bit more delicate for non-compact ones.

According to Mostow's rigidity theorem (see [Th], ch 5) volumes of hyperbolic
manifolds are homotopy invariants.  For even-dimensional ones this is
clear from the formula above.  For odd-dimensional hyperbolic
manifolds volume is
a much more interesting invariant. In this paper we will show that it is
related to
algebraic K-theory of the field of  algebraic numbers. {\it We
  will denote by $\bqq$ the subfield of all
  algebraic numbers in
$\Bbb C$}.

{\bf 2. Volumes of $(2n-1)$-dimensional hyperbolic manifolds and
  the Borel regulator on $K_{2n-1}(\bar {\Bbb Q})$}.
Let us denote by $K_*(F)$ the Quillen $K$-groups of a field
$F$.  There is the Borel regulator [Bo2]
$$
r_{m}: K_{2m-1} (\CC) \to \RR
$$

\begin{theorem} \label {Theorem 1.1}
Any $(2m-1)$-dimensional hyperbolic manifold of finite
volume $M^{2m-1}$ defines
naturally an element $\g(M^{2m-1}) \in K_{2m-1}(\bqq)\ot\QQ$ such that
 $\vo(M^{2m-1}) = r_{m} (\g(M^{2m-1}))$.
\end{theorem}

The hyperbolic volumes {vol$(M^{2n-1})$} form  a very interesting
set of numbers. If $2n-1 \geq 5$ it is discrete and, moreover,
according to
Wang's rigidity theorem [W] for any given $c \in \Bbb R$ there are only finite
number of hyperbolic manifolds of volume $\leq c$. Thank to Jorgensen and
Thurston we know that the volumes of hyperbolic 3-folds form a nondiscrete
well-odered set of ordinal type $\omega^\omega$ (see [Th]).

Theorem \ref{Theorem 1.1} together with some general conjectures in algebraic
geometry
suggests  the
highly transcendental nature of volumes of $(2n-1)$-dimensional hyperbolic
manifolds.
For example if the hyperbolic volumes $s_{1}, ...,s_{k}$ are algebraically
dependent over $\Bbb Q$ then they should be  linearly dependent
over $\Bbb Q$. Of
cource to check this is far beyond our abilities.
The first proof of theorem 1.1 is given in  chapter 3.

For an abelian group $A$ set  $A_{\Bbb Q} : = A \otimes \Bbb Q$. Set $\Bbb
Q(n) := (2\pi i)^{n}\Bbb Q$. The
complex conjugation acts on $ K_{2n-1}({ \bar {\Bbb Q}})$ and $\Bbb Q(n)$
providing the
decomposition
$$
 K_{2n-1}({ \bar {\Bbb Q}}) \otimes {\Bbb Q(n)} = (K_{2n-1}({ \bar {\Bbb Q}})
\otimes {\Bbb Q(n)})^{+}
\oplus  (K_{2n-1}({ \bar {\Bbb Q}}) \otimes {\Bbb Q(n)})^{-}
$$
In chapters 4 and 5 we will describe a different construction of an
invariant
\begin{equation} \label{din11}
\tilde \gamma(M^{2n-1})^- \in  K_{2n-1}({ \bar {\Bbb Q}}) \otimes {\Bbb
Q(n)})^{-}
\end{equation}
 which depends only on the scissor congruence class of a hyperbolic
manifold. This means that if two hyperbolic manifolds  can be cut on the
same geodesic simplices then their invariants $\tilde \gamma(M^{2n-1})^-$
coincide. We will prove that theorem 1.1 remains valid if we replace the
invariant $
\gamma(M^{2n-1})$ by $\tilde \gamma(M^{2n-1})^-$. This gives a  completely
different and more
conceptual proof this theorem.

Moreover, we will show that a formal sum of hyperbolic
polyhedrons of dimension $2n-1$ whose vertices has the
coordinates in  $\bar {\Bbb Q}$,  and which in addition
has the Dehn invariant  equal to zero,  produces an element in
$(K_{2n-1}({ \bar {\Bbb Q}}) \otimes {\Bbb
Q(n)})^{-}$. The Borel regulator on this element is equal to the
sum of volumes of the polyhedrons. In a similar situation  for
the spherical polyhedrons we get  an element in  $(K_{2n-1}({ \bar {\Bbb Q}})
\otimes {\Bbb
Q(n)})^{+}$. (The
definition of the Dehn invariant see in s. 2.2 below).
 This  was known in dimension three thanks to [DS],
[DPS], [S1], but it is
quite surprising in higher dimensions.

Theorem \ref{Theorem 1.1}  together with explicit computation
of the Borel regulator for
$K_3(\Bbb C)$ (Bl1]) and $K_5(\Bbb C)$ ([G2]) lead to
much more precise
 results about the volumes of hyperbolic 3 and
5-manifolds. To formulate them I have to
say a few words about the classical polylogarithms.

\vskip 3mm \noindent
 {\bf 3. The classical polylogarithms}.  They are defined inductively
as multivalued
analytical functions on ${\Bbb C}{\Bbb P}^1 \setminus \{0,1,\infty\}$:
$$
Li_n(z) := \int_0^z Li_{n-1}(w) \frac{dw}{w}, \qquad Li_1(z) = - \log(1-z)
$$
The function $\log z$ has a single-valued cousin $\log|z| = \re \log z$.
The coresponding function for the dilogarithm was invented by D. Wigner
and S. Bloch:
\begin{equation} \label{aa1}
{\cal L}_2(z) := \im Li_2(z) + \arg(1-z) \cdot \log|z|
\end{equation}
It is continuous on ${\Bbb C}{\Bbb P}^1$ (and, of course, real analytic on
${\Bbb C}{\Bbb P}^1
\setminus \{0,1,\infty\}$).  S. Bloch ([Bl1]) computed the
Borel regulator for $K_3(\CC)$ using ${\cal L}_2(z)$.

A single-valued version of $Li_3(z)$ that was used in [G1-G2] for an explicit
calculation of the Borel regulator for $K_5(\CC)$ looks as follows:
\begin{equation} \label{aa2}
{\cal L}_{3}(z) := \re (Li_{3}(z) - \log |z| \cdot Li_{2} (z)
+\frac13 \log^2|z|
\cdot Li_{1}(z))
\end{equation}
For arbitrary $n$ there is the following function discovered by D. Zagier [Z1]:
\begin{equation}
{\cal L}_{n}(z) := {\cal R}_{n} \left( \sum_{k=1}^n \frac{B_{k} \cdot 2^k}{k!}
Li_{n-k}(z) \cdot \log^k |z|\right)
\end{equation}

Here ${\cal R}_{n}$ means real part for $n$ odd and imaginary part for $n$ even
$B_k$ are Bernoulli numbers: $\sum_{k=0}^\infty \frac{B_k\cdot 2^k}{k!}
\cdot x^k = \frac{2x}{e^{2x}-1}$.
This function is single valued on ${\Bbb C}{\Bbb P}^1$ and
coincides with functions  (\ref{aa1}) and  (\ref{aa2})
for $n=2$ and 3.

\vskip 3mm \noindent
 {\bf 4. Volumes of hyperbolic 3 and 5-manifolds}.
Recall that the wedge  square $\La^2\CC^*$ of the multiplicative
group $\CC^*$ is defined as follows:
$$
\La^2 \CC^* := \frac{\CC^* \ot \CC^*}{\{ a \ot b + b \ot a\} } \qquad
(a,b \in \CC^*)
$$

\begin{theorem} \label {theorem 1.2}. {\it Let $M^3$ be a hyperbolic 3-manifold
of finite volume.
Then there are  numbers $z_i \in \bqq$ satisfying the condition
 \begin{equation}  \label{bloch2}
\sum (1-z_i) \land z_i = 0 \mbox { in } \La^2 \bar {\Bbb Q}^*
\end{equation}
such that}
$$
\vo(M^3) = \sum {\cal L}_2(z_i)
$$
\end{theorem}

This theorem follows from results of  Dupont--Sah [Du1], [D-S], or
Neumann-Zagier [NZ] (and was mentioned later in [Z2]).

Let $X$ be a set. Denote by $\ZZ[X]$ the free abelian group generated by
symbols $\{x\}$ where $x$ run through all elements of $X$.

Let $r(x_1, x_2, x_3, x_4) = \frac{(x_1 - x_3)(x_2 - x_4)}
{(x_1 - x_4)(x_2-x_3)}$ be the cross-ratio of 4 distinct points on a line.
Let $F$ be a field and  $R_2(F) \subset
\ZZ[P_F^1 \setminus \{0,1,\infty\}]$ be the subgroup generated by the elements
$$
\sum_{i=1}^5 (-1)^i \{r(x_1, \ldots, \hx_i,\ldots, x_5)\}\qquad
(x_i \ne x_j \in P_F^1)
$$
It is motivated by the functional
equation for the Bloch-Wigner function
$
\sum_{i=1}^5 (-1)^i{\cal L}_2 (r(z_1, \ldots, \hz_i,\ldots, z_5)) = 0,\qquad
(z_i \ne z_j \in P_\CC^1)
$. 
Moreover, any functional equation for ${\cal L}_2(z)$ is a formal
consequence of this one (see [G2], proposition 4.9).  So
$R_2(\CC)$ is the group of all functional equations
for ${\cal L}_2(z)$.
Now set
$$
B_2(F) = \frac{\ZZ[P_F^1 \setminus \{0,1,\infty\}]}{R_2(F)}
$$
Denote by $\{x\}_2$ the projection of $\{x\}$ onto $B_2(F)$.

 \begin{theorem} \label{Theorem 1.3}.  Let $M^5$ be a 5-dimensional hyperbolic
manifold
of finite volume.  Then there are algebraic numbers $z_i \in \bqq$
satisfying the condition
$$
\sum_i \{z_i\}_{2} \ot z_i = 0 \mbox{ in } B_{2}(\bqq) \ot \bqq^*
$$
such that
$$
\vo(M^5) = \sum_i {\cal L}_{3}(z_i)
$$
\end{theorem}
Theorem 1.3
(as well as  theorem 1.2) follows from theorem 1.1 and the main results of
[G2].

{\bf 5. Volumes of hyperbolic (2n-1)-manifolds}. One can define for
arbitrary $n$ a certain subgroup ${\cal R}_{n}(F) \subset
\ZZ[P_F^1]$ which for $F=\CC$ is the subgroup of {\it all\/} functional
equations for the $n$-logarithm ${\cal L}_{n}(z)$ (see s. 1.4 of [G2]).  Set
$$
{\cal B}_{n}(F) := \frac{\ZZ[P_F^1]}{{\cal R}_{n}(F)}
$$
Let us define  the  homomorphism
$$
\de_n: \Bbb Z[P^1_F] \to {\cal B}_{n-1}(F)\ot F^*
\qquad \{x\} \mapsto \{x\}_{n-1}\ot x
$$
One can show that $\de_n ({\cal R}_n(F) )=0$ (see [G2]), so we
get a homomorphism
$$
\de_n: {\cal B}_{n}(F) \to {\cal B}_{n-1}(F)\ot F^*
$$
\begin{conjecture} \label{Conjecture 1.4} {\it Let $M^{2n-1}$ be an
$(2n-1)$-dimensional hyperbolic
manifold of finite volume.  Then there are algebraic numbers $z_i \in \bqq$
 satisfying the condition $(n \ge 3)$
$$
\sum_i \{ z_i\}_{n-1} \ot z_i = 0 \mbox{ in } {\cal B}_{n-1}(\bqq) \ot
\bqq^*
$$
such that}
$$
\vo(M^{2n-1}) = \sum_i {\cal L}_{n}(z_i)
$$
\end{conjecture}
This conjecture   is a consequence of a version of Zagier's conjecture
[Z1] and  theorem 1.1.

\vskip 3mm \noindent
{\bf 6. The Chern-Simons invariants of compact hyperbolic
  (2n-1)-manifolds}. The Chern-Simons invariant $CS(M^{2n-1})$ takes
the values in $S^{1} = {\Bbb R}/{\Bbb Z}$.
There is the Beilinson's regulator
\begin{equation} \label{44}
r^{Be}_{n}: K_{2n-1}(\bar {\Bbb Q}) \longrightarrow
\Bbb C/(2\pi i)^{n} \Bbb Q
\end{equation}
Therefore $(2\pi)^{n}  r^{Be}_{n}$ provides homomorphisms
$$
(K_{2n-1}({ \bar {\Bbb Q}})\otimes \Bbb Q(n))^{-} \longrightarrow \Bbb R
\qquad (K_{2n-1}({ \bar {\Bbb Q}})\otimes \Bbb Q(n))^{+}
\longrightarrow \Bbb R/ \Bbb Q
$$
For an element $\gamma \in K_{2n-1}({ \bar {\Bbb Q}})$ let
$\gamma^{\pm}$ be its components in
$(K_{2n-1}({ \bar {\Bbb Q}})\otimes \Bbb Q(n))^{\pm}$

Suppose that $M^{2n-1}$ is a compact
  hyperbolic manifold.
Then one can show that
$$
CS(M^{2n-1}) = c_{n} r^{Be}_{n}(\gamma(M^{2n-1})^{+})
$$
where $c_n$ is a nonzero universal constant. For 3-manifolds this
was known: see [D2] if the manifold is compact, and
[Y], [N], [NJ] if it has cusps. The proof of this result will
be given
elsewhere.
Notice that according to our conventions vol$ (M^{2n-1}) =
r^{Be}_{n}(\gamma(M^{2n-1})^{-})$.

As far as I know the Chern-Simons invariants of noncompact hyperbolic manifolds
 was defined only for 3-manifolds ([Me]). So one can use the formula above
 as a {\it definition} of the $CS(M^{2n-1})$ for
 noncompact hyperbolic manifolds.

{\bf 7. The invariants  ${\tilde \g}(M^{2n-1})$ and  ${\tilde
    \g}(M^{2n-1})^-$}. There are
two completely different points
of view on algebraic K-theory.

i). In Quillen's
definition $K_{n}(F)_{\Bbb Q}$   is the quotient of
$H_{n}(GL(F),\Bbb Q)$ modulo the subspace
of decomposable elements. This is how elements $\g(M^{2m-1})$
constructed in chapter 3. The construction is straightforward
for
3-dimensional compact hyperbolic manifolds: the fundamental class
of such a manifold provides an element in $H_3(BSO(3,1), \Bbb
Z)$; using the local isomorphism between $SO(3,1)$ and $SL(2,\Bbb
C)$ we get a clas in $H_3(SL(2,\Bbb
C)^{\delta})$ and thus an element in $H_3(SL(\Bbb
C)^{\delta})$ via the embedding $SL_2 \hookrightarrow SL$. It is
interesting that for compact hyperbolic manifolds of dimension
$\geq 5$ one have to use the halfspinor representation of
$SO(2n-1,1)$ in order
to get,  starting from the fundamental
class of $M^{2n-1}$, an interesting  class in $H_{2n-1}(GL(\Bbb C))$. (It seems
that the other
fundamental representations of $SO(2n-1,1)$ lead to zero classes in
$K_{2n-1}(\Bbb C)$). For the noncompact hyperbolic manifolds
$H_{2n-1}(M^{2n-1}) = 0$, and the construction of the invariant
$\gamma(M^{2n-1})$ becomes rather delicate,
see chapter 3.

ii). Let $gr^{\gamma}_{n}K_{m}({{\Bbb C}})$ be the
graded quotients of the
 $\gamma$-filtration on $K_{m}({\Bbb C})$.
 According to A.A.Beilinson elements of
$gr^{\gamma}_{n}K_{2n-1}({\Bbb C})$ should have an interpretation as motivic
extensions of ${\Bbb Q}(0)$ by ${\Bbb Q}(n)$ (here ${\Bbb Q}(n)$ is the $n$-th
tensor power of the Tate motive ${\Bbb Q}(1)$; the latter is the inverse
 to the motive $H^{2}({\Bbb C}{\Bbb P}^{1})$ ).

There is a natural map $ K_{2n-1}({ \bar {\Bbb Q}})
 \rightarrow
gr^{\gamma}_{n}K_{2n-1}({\Bbb C})$. It is expected to be an
isomorphism modulo torsion (the rigidity conjecture).

In the chapter 4 (see also chapter 2 for an introduction) for any
hyperbolic manifold of finite
volume $M^{2n-1}$
I will give a simple geometrical construction of a mixed Tate motif
$m(M^{2n-1})$ that fits in the following exact sequence
(in the
category of mixed Tate motives):
$$
0 \longrightarrow {\Bbb Q}(n)  \longrightarrow m(M^{2n-1})
\longrightarrow {\Bbb Q}(0) \longrightarrow 0
$$
Namely,
I will construct an  extension in the category of mixed Hodge structures that
will be clearly of algebraic-geometric origin.

Moreover, it will be clear from the construction that the $\Bbb R$-part  of the
regulator map on the class of $Ext^1$ defined by
this extension coincides with the volume of $M^{2n-1}$.

This extension  represents an element
$$
{\tilde \gamma}(M^{2n-1})^{-} \in (K_{2n-1}(\bar
{\Bbb Q})\otimes \Bbb Q(n))^-$$
One should have
${\tilde \gamma}(M^{2n-1})^{-} = \gamma(M^{2n-1})^{-}$.

The element ${\tilde \gamma}(M^{2n-1})^{-}$ depends only on the scissor
congruence class of the manifold $M^{2n-1}$ just by the construction. Therefore
the invariant $\gamma(M^{2n-1})^{-}$ is responsible for the  scissor
congruence class of the hyperbolic manifold $M^{2n-1}$ (see
s.2.1-2.2 below). Moreover, the following
three conditions {\it should be} equivalent:

a) The volumes of hyperbolic manifolds $M_{1}^{2n-1}$ and $M_{2}^{2n-1}$ are
coincide.

b)$\gamma(M_{1}^{2n-1})^{-}  =  \gamma(M_{2}^{2n-1})^{-}$.

c)$M_{1}^{2n-1}$ and $M_{2}^{2n-1}$ are scissor congruent.

{\bf 8. The structure of the paper}. In chapter 2 we define Dehn
complexes, which generalize the  scissor congruence
groups and the classical Dehn invariant of a
polyhedron, in the hyperbolic, spherical and euclidean
geometries. We formulate  conjectures relating the cohomology of these
complexes to algebraic K-theory of $\Bbb C$ and (in the euclidean
case)
 to Kahler differentials $\Omega^i_{\Bbb R/\Bbb Q}$. A part of these
 conjectures for noneuclidean polyhedrons whose vertices has
 coordinates in a number field will be
 proved in chapters 4 and 5. The crucial result is an
 interpretation of the Dehn invariant on the language of mixed
 Hodge structures given in chapter 4.

The chapter 2 can be considered as  a continuation of the
introduction: it contains some important definitions and
formulations of the results and conjectures as well as  ideas of
some
of constructions, but we avoid any technical discussion.

More generally, in chapter 4
we will construct some commutative graded Hopf algebra
$S(\Bbb C)_{\bullet}$ which  is  dual to
the universal enveloping algebra of a certain remarkable
negatively graded (pro)-Lie algebra. I believe that
the category of finite dimensional graded modules over this Lie
algebra  is equivalent to the category of mixed Tate motives
over $\Bbb C$ (compare with the Hopf algebra considered in
[BGSV]).
In particulary this means that  the cohomology of $S(\Bbb
C)_{\bullet}$ should coincide with appropriate pieces of
the algebraic K-theory of ${\Bbb C}$.

Chapter 5 starts from certain known to experts technical results
about triangulated
categories and their abelian hearts. We need them to finish the
proof of theorem \ref{Theorem 2.4}.

\section {Dehn complexes in classical geometries and their homology}

\vskip 3mm \noindent
{\bf 1. The scissor congruence groups }.
There are three classical geometry: euclidean, hyperbolic and
spherical. In each of them one can consider geodesic simplices and ask
the following question

{\bf Hilbert's third problem}. {\it Suppose two geodesic simplices has the same
volume. Is it possible to cut one of them on geodesic simplices
and putting them together in a different way get another simplex?}

To formulate the problem in a more accurate way we need to introduce the
scissor congruence groups.

Let $V^n$ be an $n$-dimensional space with one of the classical
geometries,
 i.e. $V$ is the hyperbolic space ${\cal H}^{n}$, spherical space $S^n$
or Euclidian space $E^n$

Any $n + 1$
points $x_{0}, ... , x_{n}$ in the space $V^{n}$ define a
geodesic simplex $I(x_{0}, ... , x_{n})$ with
vertices in these points. Let us denote by ${\cal P}( V^{n})$ the
abelian group generated by symbols $\{I(x_{0}, ... , x_{n}), \alpha \}$ where
$\alpha$ is an orientation of $V^{n}$, subject to the following
relations:

a) $\{I(x_{0}, ... , x_{n}), \alpha \} = 0$ if $x_{0}, ... , x_{n}$ lie in a
geodesic hyperplane.

b)$\{I(x_{0}, ... , x_{n}), \alpha \} = \{I(gx_{0}, ... , gx_{n}), g\alpha \}$
for any element $g$ form the group of automorphisms of the corresponding
geometry.

c)$\{I(x_{\sigma(0)}, ... , x_{\sigma(n)}), \beta \} =
-(-1)^{|\sigma|}
\{I(x_{0}, ... , x_{n}), \alpha \}$ where $\sigma$ is a permutation and $\beta$
is another orientation of $V^{n}$.

d) $\sum _{i = 0}^{ n+1} (-1)^{i}\{I(x_{0}, ... ,\hat x_{i}, ... , x_{n+1}),
\alpha \}
= 0$ for any $n$+2 points $x_{i}$ in $V^{n}$.

The automorphism groups are: the group $O(n,1)$ for the
hyperbolic case, $O(n+1)$ for the spherical geometry and the semidirect
product of $O(n)$ and translation group $E^n$ for the Euclidian geometry

\begin{lemma} \label{Lemma 2.1}
${\cal P} ( S^{2n}) = 0$.
\end{lemma}

{\bf Proof}. Points $-x_0,x_0,x_1,...,x_{2n-1}$ belong to a geodesic hyperplane
in $S^{2n}$. Therefore $\{I(-x_0,x_0,x_1,...,x_{2n-1}),\alpha\} =0$.

Using this fact we see that the  additivity axiom d) for $(-x_{0},x_{0},,
x_{1}, ... , x_{2n})$ together with a)
implies that
$(I(x_{0},x_{1}, ... ,x_{2n}), \alpha) =
(I(-x_{0},x_{1}, ... ,x_{2n}), \alpha)$
and hence $(I(x_{0},x_{1}, ... ,x_{2n}), \alpha) = (I(-x_{0},-x_{1}, ...
,-x_{2n}), \alpha)$. The relation b) for the antipodal involution is
$(I(x_{0},x_{1}, ... ,x_{2n}), \alpha) = (I(-x_{0},-x_{1}, ... ,-x_{2n}),
\beta)$. Applying c) we get lemma 2.1.

The volume of a geodesic simplex provide  homomorphisms
$$
vol_{{\cal H}}: {\cal P} ( {\cal H}^{n}) \rightarrow \Bbb R,
\qquad vol_{S}: {\cal P}
( S^{2n-1}) \rightarrow \Bbb R/ \Bbb Z,    \qquad  vol_{E}: {\cal P} ( E^{n})
\rightarrow \Bbb R
$$

 It is not hard to see (and in the euclidean case was known to Euclid)
that the volume  provides isomorphisms
$$
{\cal P} ( S^{1}) = \Bbb R/ \Bbb Z , \qquad {\cal P}
( {\cal H}^{1}) = {\cal P} ( {\cal H}^{2}) = \Bbb R, \qquad {\cal P} (
E^{1}) = {\cal P} ( E^{2}) = \Bbb R
$$

A more accurate formulation of the Hilbert's third problem is the following
question: is it true that the volume homomorphism is injective in dimension
three?

The negative answer to this problem was given by Max Dehn (1898). He
discovered that in dimension
$\geq 3$ the volume
{\it does not} separate the elements of the scissor congruence group,
because a
new phenomena appears:

{\bf 2. The Dehn  invariant}.   This is a homomorphism
\begin{equation}  \label{dehnhh}
D_{n}^{ V}: {\cal P} (  V^{n}) \longrightarrow \bigoplus ^{n-2}_{i =1} {\cal P}
(  V^{i}) \otimes {\cal P} ( S^{n-i-1})
\end{equation}
Here, as usual, $V^{n}$ can be  ${\cal H}^n$, $S^{n}$ or $E^{n}$.

It is defined as follows. Let $({I, \alpha})$
be a generator of ${\cal P}( V^{n})$. For each $i$-dimensional edge $A$
of a (geodesic) simplex $I$ consider the corresponding geodesic $i$-plane.
This plane inherits the same type of geometry as $V^n$ (i.e. it is a
hyperbolic plane in hyperbolic geometry etc.) and vertices of the edge $A$
define a geodesic
simplex in it.
The numeration of its vertices is induced by the one of $I$. Choose an
orientation $\alpha_{A}$ of this $i$-plane. We get an element of ${\cal P}
( V^{i})$  that will be denoted as $I_{A}$.

Now look at a $(n-i)$-plane $A'$ orthogonal to $A$ and intersecting it at a
certain point $e $. This is  a Euclidian plane in all
geometries. The intersection of the
sphere in $A'$ centered in $e$ with $i+1$-dimensional edges of $I$ containing
$A$ defines a spherical simplex. Its vertices has a natural numeration induced
by the numeration of the vertices if $I$. Finally, there is an orientation
$\alpha_{A'}$ of $A'$ such that the orientation of ${\cal H}^{n}$ defined by
$\alpha_{A}$ and
$ \alpha_{A'}$ coincides with the orientation $\alpha$. So we get an element
$I_{A'} \in
{\cal P} ( S^{n-i})$. By definition
$$
D^{H}_{n} (I,\alpha):= \sum_{A} I_{A} \otimes I_{A'}
$$
where the sum is over all edges of the simplex $I$ of dimension $i$,
$0<i<n$

Useful general references on this subject are [Sah 1] and [C]. In
particular according to S.H.Sah $\oplus {\cal P} ( S^{n})$ is a Hopf algebra
with spherical
Dehn invariant $D^{S}_{n}$   as a comultiplication; hyperbolic
(respectively euclidean) Dehn invariant  (\ref{dehnhh})
provides
$ {\cal P} ( {\cal H}^{n})$ (respectively $ {\cal P} ( E^{n})$ ) with a
structure of a comudule over
it. Notice that our groups ${\cal P} ( S^{n})$ are different (smaller) from the
ones defined in [Sah]. For those groups , for example, one has ${\cal P}
( S^{2n}) ={\cal P} ( S^{2n-1})$ .

{\bf 3. Scissor congrunce class of a hyperbolic manifold}. Any
hyperbolic manifold of finite volume can be cut into a finite number of
 geodesic simplices (see, for example, chapter 2 below). If the manifold is
compact
these simplices have vertices in ${\cal H}^{n}$. Therefore one can consider the
sum of corresponding elements in ${\cal P}( {\cal H}^{n})$. If the manifold is
noncompact, the vertices of these simplices might be in the absolute
$\partial {\cal H}^{n}$. However one can introduce the scissor congruent group
${\cal P}({\bar {\cal  H}}^{n})$ generated by simplices with vertices at
${\bar {\cal  H}}^{n} = {\cal  H}^{n} \cup  \partial {\cal H}^{n}$.
It turns out that

\begin{proposition} \label{Proposition 2.2}
([Sah 2],    ) The natural inclusion
\begin{equation} \label{sah}
{\cal P}( {\cal H}^{n}) \hookrightarrow {\cal P}({\bar {\cal  H}}^{n})
\end{equation}
is an isomorphism.
\end{proposition}

The sum of the elements of
${\cal P}( {\cal H}^{n})$ corresponding to these simplices clearly does not
depend on the cutting. Therefore {\it any} hyperbolic manifold
produces an element
$s(V^{n}) \in {\cal P}( {\cal H}^{n})$

\begin{proposition} \label{Proposition 2.3}
The Dehn invariant of
   $s(M^{n})$ is equal to zero.
\end{proposition}

{\bf Proof}. This is quite clear for compact
manifolds. Namely,consider a $k$-dimensional edge $A$ of the
triangulation of $V^{n}$ on geodesic simplices. Let $I^{i}$ be
the set of all simplices containing the edge $A$. Each simplex
$I^i$ defines an element $I_{A'}^{i} \in {\cal P} ( S^{n-k})$:
the ``inner angle''
at the edge $A$ (see above). It follows from the very definition that
 $I_{A}$ will appear in formula for $D_{n}^{H}$ with factor
$\sum_{i} I_{A'}^{i}$. But $\sum_{i} I_{A'}^{i} = [S^{n-k}] = 0$ in
${\cal P} ( S^{n-k})$ because $M^{n}$ is a manifold without boundary.

For noncompact manifolds proposition 2.3
 is less obvious and we have to proceed as follows. One can
 define the (extended) Dehn invariant
\begin{equation} \label{s}
D_{2n-1}^{ H}: {\cal P} (\bar {\cal H}^{2n-1}) \longrightarrow
\bigoplus ^{n-1}_{i =1} {\cal P}
(\bar {\cal H}^{2i-1}) \otimes {\cal P} ( S^{2(n-i)-1})
\end{equation}
providing that the homomorphism  $(\ref{sah})$  commutes with the Dehn
invariant. Namely, the definition of the components of $(\ref{s})$
corresponding to edges of dimension bigger then 1 is {\it verbatim} the same.
For 1-dimensional edges let us use the following Thurston's regularization
procedure. Let $I(x_{1}, ..., x_{2n})$ be a simplex in ${\cal P} (\bar {\cal
H}^{2n-1})$   .For each infinite
vertex $x_{i}$  delete a small horoball centered at
$x_{i}$. Denote by $x_{i}(j)$ the intersection point of the edge $x_{i} x_{j}$
with the horosphere. If $x_{k}$ is a point inside hyperbolic space, set
$x_{k}(j) = x_{k}$. Then the component of $(\ref{s})$ corresponding to
the edge $x_{i}
x_{j}$ is the tensor product of the 1-dimensional hyperbolic simplex
$x_{i}(j) x_{j}(i)$ and the spherical $2n-1$ simplex that appears in the usual
definition of Dehn invariant (`` inner angle'' at the edge $x_{i}x_{j}$,
see above). The expression we get
does not depend on the choice of horoballs because of the following
reason. Consider a simplex in Euclidian space ${\Bbb R}^{2n-2}$.
Each its vertex $v$ defines a spherical simplex and hence an element of
${\cal P}(S^{2n-2})$ (``inner angle'' at $v$). Then the  sum of
the elements corresponding to all vertices of the simplex is zero (this
generalises  the fact that the sum of angles of an euclidean triangle is
$\pi$ and true only for simplices in even-dimensional spaces).

After this the proof for noncompact hyperbolic manifolds
is the same as for compact ones.

{\bf 4. Kernel of the Dehn invariant in non-euclidean geometry and
$K_{2n-1}({\bar {\Bbb Q}})\otimes \Bbb Q$ }.
\begin{conjecture} \label{Conjecture 2.4}
There are canonical injective  homomorphisms
$$
KerD^{H}_{2n-1} \otimes {\Bbb Q} \hookrightarrow (K_{2n-1}({ \bar {\Bbb Q}})
\otimes \Bbb Q(n))^{-}
$$
$$
KerD^{S}_{2n-1} \otimes {\Bbb Q} \hookrightarrow (K_{2n-1}({ \bar {\Bbb Q}})
\otimes \Bbb Q(n))^{+}
$$
such that the following diagrams are commutative:
\begin{center}
\begin{picture}(100,70)
\put(-50,50){$
KerD^{H}_{2n-1} \otimes {\Bbb Q}
$}
\put(33,53){\vector(1,0){80}}
\put(120,50){$
(K_{2n-1}({ \bar {\Bbb Q}}) \otimes \Bbb Q(n))^{-},
$}
\put(55,-15){$
\Bbb R
$}
\put(55,55){$ $}
\put(-10,40){\vector(1,-1){45}}
\put(130,40){\vector(-1,-1){45}}
\put(110,10){$r_{n}$}
\put(-5,10){$vol$}
\end{picture}
\end{center}
\vskip 4mm

\begin{center}
\begin{picture}(100,70)
\put(-50,50){$
KerD^{S}_{2n-1} \otimes {\Bbb Q}
$}
\put(33,52){\vector(1,0){80}}
\put(120,50){$
(K_{2n-1}({ \bar {\Bbb Q}}) \otimes \Bbb Q(n))^{+},
$}
\put(49,-17){$
{\Bbb R}/{\Bbb Z}
$}
\put(55,55){$ $}
\put(-10,40){\vector(1,-1){45}}
\put(130,40){\vector(-1,-1){45}}
\put(110,10){$r^{Be}_{n}$}
\put(-5,10){$vol$}
\end{picture}
\end{center}
\end{conjecture}

\vskip 5mm
Let ${\cal P}({\cal H}^{2n-1},{\bar {\Bbb Q}})$ (respectively ${\cal P}
(S^{(2n-1)},{\bar {\Bbb Q}})$ be the subgroup of the
hyperbolic (resp. spherical) scissor congruence group generated by the
simplices
defined over ${\bar {\Bbb Q}}$ (i.e. their vertices have
coordinates in ${\bar {\Bbb Q}}$.

\begin{theorem} \label{Theorem 2.4}
 {\it There are the following commutative diagrams }

\begin{center}
\begin{picture}(100,70)
\put(-90,50){$
KerD^{H}_{2n-1}|_{{\cal P}({\cal H}^{2n-1},{\bar {\Bbb Q}})}\otimes {\Bbb Q}
$}
\put(33,53){\vector(1,0){80}}
\put(120,50){$
(K_{2n-1}({ \bar {\Bbb Q}}) \otimes \Bbb Q(n))^{-},
$}
\put(55,-15){$
\Bbb R
$}
\put(55,55){$ $}
\put(-10,40){\vector(1,-1){45}}
\put(130,40){\vector(-1,-1){45}}
\put(110,10){$r_{n}$}
\put(-5,10){$vol$}
\end{picture}
\end{center}
\vskip 4mm

\begin{center}
\begin{picture}(100,70)
\put(-90,50){$
KerD^{S}_{2n-1} |_{{\cal P}(S^{2n-1}),{\bar {\Bbb Q}} )}  \otimes {\Bbb Q}
$}
\put(33,52){\vector(1,0){80}}
\put(120,50){$
(K_{2n-1}({ \bar {\Bbb Q}}) \otimes \Bbb Q(n))^{+},
$}
\put(49,-17){$
{\Bbb R}/{\Bbb Z}
$}
\put(55,55){$ $}
\put(-10,40){\vector(1,-1){45}}
\put(130,40){\vector(-1,-1){45}}
\put(110,10){$r^{Be}_{n}$}
\put(-5,10){$vol$}
\end{picture}
\end{center}
\end{theorem}
\vskip 5mm
The proof of this theorem is based on the following  ideas:

1. A $(2n-1)$-dimensional non-euclidean geodesic simplex $M$ defines a mixed
Tate motive $\gamma(M)$
with an additional data: $n$-framing. See s. 2.6 below and all
the details in chapter 4.

2. The Dehn invariant also has a motivic interpretation and moreover
each element $[M] \in KerD^{H}_{2n-1}|_{{\cal P}({\cal H}^{2n-1},{\bar {\Bbb
Q}})}\otimes {\Bbb Q}$ defines an exact sequence in the category of
mixed Tate motives:
$$
0 \longrightarrow {\Bbb Q(n)}_{{\cal }} \longrightarrow
\gamma([M]) \longrightarrow
{\Bbb Q(0)}_{{\cal }}\longrightarrow 0
$$
i.e. an element of $Ext^{1}_{\cal M}(\Bbb Q(0)_{\cal M},\Bbb Q(n)_{\cal M})$.
See chapter 4.

3. We already have an abelian category of mixed Tate motives over
a number field. This follows from the existence of a certain
triangulated (not abelian!) category of mixed motives over a
field $F$ constructed by Levine [L] and Voevodsky [V] in which
the Ext-groups  isomorphic to appropriate parts of the K-theory,
and some formal arguments that use the Borel theorem, see the details in
chapter 5.

The following conjecture, I guess, express a general expectation about the
kernel of euclidean Dehn invariant. It tells us that the Dehn invariant and
the volume homomorphism separate all elements of the euclidean scissor
congruence group. Some speculations about its generalization  and an
``explanation''   see below in s. 2.7.
\begin{conjecture} \label{eucl}
$Ker D^E = \Bbb R$
\end{conjecture}

{\bf 5. Dehn complexes in non euclidean geometry and algebraic K-theory of
$\Bbb C$}.
The next natural question is what can we say about the cokernel of the
Dehn invariant. First of all it turns out that the Dehn invariant in
higher (bigger then 3) dimensions is only the beginning of a certain
complex, which I will call the Dehn complex. Namely,
 consider the following complexes. (Here, again, $V$ means one of
the three geometries)
\begin{eqnarray} \label{dehn1}
{\cal P}_V^{\bullet}(n) \qquad :{\cal P}( V^{2n-1}) \stackrel {d_{v}}
{\longrightarrow}
\oplus_{i_{1} + i_{2} = 2n-1} {\cal P}( V^{i_{1}})
\otimes {\cal P}(S^{i_{2}}-1)
\stackrel {d_{v}}
{\longrightarrow} ...   \nonumber\\
\stackrel {d_{v}}
{\longrightarrow } \oplus_{i_{1} + ... +i_{k}= 2n-1}
{\cal P}( V^{i_{1}}) \otimes {\cal P}(S^{i_{2}-1})\otimes ... \otimes
{\cal P}(S^{i_{k}-1})
\stackrel {d_{v}}
{\longrightarrow} ...
\end{eqnarray}

Here the first group is placed in degree 1 and the differentials
have degree $+1$;
\begin{equation}
d_{v} = D^{V} \otimes id \otimes ... \otimes id -  id \otimes
D^{S} \otimes id \otimes ... \otimes id + ... \pm
\otimes id \otimes ... \otimes id \otimes D^{S}
\end{equation}

{\bf Remark}. The spherical Dehn invariant provides $\oplus {\cal
  P} (S^{2j-1})$ with a
structure of a coalgebra.
Further, the Dehn invariant provides $\oplus{\cal P}(V^{2i-1})$ with a
structure of a  comodules over this
coalgebra.
The complexes
above are just the cobar complexes
computing the degree $n$ part of the cohomology of the
coalgebra $\oplus {\cal
  P} (S^{2n-1})$ with coefficients in comodule
$\oplus{\cal P}({\cal H}^{2n-1})$.

\begin{conjecture} \label{conjecture 2.5}. {\it There are canonical
homomorphisms}
\begin{equation} \label{1111}
H^{i}({\cal P}_{\cal H}^{\bullet}(n)) \longrightarrow
(gr^{\gamma}_{n}K_{2n-i}(\Bbb C)\otimes \Bbb Q(n))^{-}
\end{equation}
\begin{equation} \label{11111}
H^{i}({\cal P}_{S}^{\bullet}(n)) \longrightarrow (gr^{\gamma}_{n}K_{2n-i}(\Bbb
C)\otimes \Bbb Q(n))^{+}
\end{equation}
\end{conjecture}

The following conjecture is due to D. Ramakrishnan and generalizes Milnor's
conjecture about the values of the Lobachevsky function ( see 7.1.2 in [R])

\begin{conjecture} \label{Conjecture 2.6}
The Beilinson regulator is injective modulo
torsion.
\end{conjecture}
Conjectures  (\ref{Conjecture 2.4}), (\ref{Conjecture 2.6}) would imply that in
hyperbolic and spherical geometry
 the

\vskip 3mm \noindent
{\bf The Extended Hilbert Third Problem}. {\it Do the Dehn invariant and the
volume separate all points of scissor congruence groups?}

should have an affirmative answer.

\vskip 3mm \noindent
{\bf Problem} {\it  Is it true that the homomorphisms (\ref{1111}),
  (\ref{11111}) are isomorphisms modulo torsion?}

For $n=2$ the answer is yes. This follows from the
 results of J.Dupont and S.H.Sah(see [D], [DS1], [DPS], [Sah3]).Unfortunately
their methods use essentially  classical isomorphisms between simple Lie groups
in low dimensions (like
the local isomorphism between $SO(3,1)_{0}$ and $SL_{2}(\Bbb
 C)$) and a lot of
arguments ``ad hoc''
that one does not see how to generalize to higher dimensions,
even hypothetically.
Their approach deal with the homology of orthogonal groups and  the relation
with algebraic $K$-theory was
quite mysterious.

In the next section and in the chapter 4 below I will give a
motivic interpretation
of Dehn complexes. After this algebraic $K$-theory shows up immediately thanks
to Beilinson's conjectures about the category of mixed Tate
motives. In particularly this approach
clarifies why the dilogarithm appears in the classical computations of volumes
of geodesic 3-simplices (see, for example, [Co] and [M2]).

\vskip 3mm \noindent
{\bf 6. Geodesic simplices in non-Euclidian geometry and mixed Tate motives}.
In the chapters 4 and 5 I will deduce conjectures (\ref{Conjecture 2.4}) and
(\ref{conjecture 2.5})
from standard conjectures about mixed Tate motives. The key idea
is that any geodesic simplex $M$ in spherical or hyperbolic
geometry define
certain mixed  Tate motive. To get its Hodge realization in the case , say,
hyperbolic geometry one has to proceed as follows.

In the Klein model the Lobachevsky space ${\cal H}^m$ is realized
as the interior of a ball in $\Bbb R^m$ and the distance $\rho(P_1, P_2)$ is
defined as
$|\log r(Q_1, Q_2, P_1, P_2)|$ where $Q_1$, $Q_2$ are the intersection points
of the line $P_1P_2$ with the absolute: the sphere $Q$, and $Q_1,
P_1, P_2, Q_2$ is the order of the points on the line.

Then geodesics are straight lines and so a geodesic simplices are
just the usual
ones inside $Q$.

Let us complexify and compactify this picture.  We will get
${\Bbb C}{\Bbb P}^m$ together
with a quadric $Q$ corresponding to the absolute and a collection of
hyperplanes $M=(M_1, \ldots, M_{m+1})$ corresponding to faces of
a geodesic simplex.

For any nondegenerate quadric $Q \subset {\Bbb C}{\Bbb P}^m$ and
a simplex $M=(M_1,
\ldots, M_{m+1})$ in generic position with respect to $Q$
let us denote by $H^m({\Bbb C}{\Bbb P}^m \setminus Q, M)$ the
$m$-th cohomology group
of ${\Bbb C}{\Bbb P}^m\setminus Q, \Bbb Q)$ modulo $Q\setminus
(Q\cap M)$.   This space has
canonical mixed Hodge structure, which we denote as $h(Q,M)$.
The weights are $0,2,\ldots, 2\cdot \left[ \frac{m-1}{2}\right]$
and the corresponding graded quotients of weight
$2j$ are isomorphic to direct sum of the Tate structures
$\QQ(-j)$. Such mixed Hodge structures are called mixed
Hodge-Tate structures. In particularly a geodesic simplex
produces a Hodge-Tate structure.

In complete analogy with this a quadric $Q \in {\Bbb C}{\Bbb
  P}^{m}$ without points in  ${\bf RP}^{m}$  defines an elliptic
geometry in
${\bf RP}^{m}$; the geodesic simplices are just the usual ones.

Now let us suppose that $m=2n-1$.
There is a canonical, up to a sign, meromorphic $2n-1$ - form $\omega_{Q}$ in
${\Bbb C}{\Bbb P}^{2n-1} \backslash Q$ with polar singularity
along $Q$. If we choose coordinates
$x_{1}, x_{2}, ... , x_{2n}$ in
$\CC^{2n}$ and an equation $\tq = \sum_{i,j} q_{ij} x_i x_j = 0$
of the quadric $Q$,  then
\begin{equation}
w_Q := 1 /{(2\pi i)^n} \cdot \sqrt{\det \tq} \cdot\frac{\sum_{i=1}^{2n} (-1)^i
x_i dx_i
\land \ldots \land \wh{dx_i} \land \ldots \land dx_{2n}}{\tq^n}
\end{equation}
It does not depend on the choice of the coordinates $x_i$ and the equation
$\tq$.

Let $\Delta_M$ be a relative cycle representing a generator
of $H_{2n-1}({\Bbb C}{\Bbb P}^{2n-1}, M)$. Set
\begin{equation}
v(Q,M) = \int_{\Delta_M} \omega_{Q}
\end{equation}

If a pair $(Q,M)$ corresponds to a geodesic simplex $S$ in ${\cal H}^{2n-1}$
then there is a natural choice for the relative cycle $\Delta_M$ and
$v(Q,M) = \vo(S)$. So the volume of an odd-dimensional geodesic simplex is the
maximal period of the corresponding mixed Hodge structure.

Mixed Tate motives over ${\Bbb C}$ should be objects of certain abelian  $\Bbb
Q$-category $\cal M$. There should be a realization functor from
$\cal M$ to the category of mixed Hodge-Tate structures ${\cal H}_{T}$.
According to general philosophy mixed Tate motives
can be considered as those of mixed Hodge-Tate structures that can be realised
in cohomology of algebraic varieties. Therefore the
 mixed Hodge-Tate structure $h(Q,M)$ constructed above definitely corresponds
to a certain mixed Tate motive.

However Ext-groups between even simplest mixed Tate motives
$\Bbb Q(n)_{\cal M}$ in the category ${\cal M}$ should be quite different from
those of their Hodge counterparts $\Bbb Q(n)_{{\cal H}_{T}}$. For example,
 one should have ([B1])
\begin{equation} \label{ext111}
Ext^{1}_{\cal M}(\Bbb Q(0)_{\cal M},\Bbb Q(n)_{\cal M}) \otimes \Bbb Q
= gr^{\gamma}_{n}K_{2n-1}(\Bbb C) \otimes \Bbb Q
\end{equation}
while in the category of mixed Hodge-Tate structures one have
\begin{equation} \label{br11}
Ext^{1}_{{\cal H}_{T}}(\Bbb Q(0)_{{\cal H}_{T}},\Bbb Q(n)_{{\cal H}_{T}}) =
\Bbb C/(2\pi i)^{n} \Bbb Q
\end{equation}
and the realization functor provides a homomorphism
\begin{equation}
Ext^{1}_{\cal M}(\Bbb Q(0)_{\cal M},\Bbb Q(n)_{\cal M}) \longrightarrow
Ext^{1}_{{\cal H}_{T}}(\Bbb Q(0)_{{\cal H}_{T}},\Bbb Q(n)_{{\cal H}_{T}})
\end{equation}
that should coincide with the Beilinson's regulator (\ref{44})
after identifications (\ref{ext111}) (\ref{br11}).

The mixed Hodge-Tate structure $h(Q,M)$ has the following
additional data ($n$ -framing): the form $\omega_{Q}$ defines a vector
$[w_Q] \in gr^W_{2n}h(Q,M)$ and a generator of $H_{2n-1}(\CC P^{2n-1},M)$
can be viewed as a functional in $gr^W_0 h(Q,M)$. (Compare with [BGSV]).

$n$-framed mixed Hodge-Tate structures
form an abelian group ${\cal H}_{n}$ (see [BGSV], or s.4.2). There is a
homomorphism of a
scissor congruence group in a $(2n-1)$- dimensional non-Euclidian space to
${\cal H}_{n}$ defined on generators as described
above.

The next important idea is that the Dehn invariant has a natural
interpretation in the language of mixed Hodge-Tate structures (theorem 4.8).
Moreover, any element in the kernel of the Dehn invariant
corresponds to an element of
$Ext^{1}{{\cal H}_{T}}(\Bbb Q(0)_{{\cal H}_{T}},\Bbb Q(n)_{{\cal
    H}_{T}})$. It is clearly of motivic origin and hence gives us
an element of (\ref{ext111}). Therefore we get a map predicted by
conjecture 2.4 in the hyperbolic case.

{\bf 7. Some speculations about euclidean Dehn complexes and  mixed Tate
motives over the dual
numbers}.
Let $\Omega^i_{\Bbb R/\Bbb Q}$ be the Kahler differential i-forms of the
field $\Bbb R$
\begin{conjecture} \label{deuc}
$H^i({{\cal P}}^{\bullet}_{E}(2n-1))\otimes \Bbb Q = \Omega^{i-1}_{\Bbb R/\Bbb
Q}$
\end{conjecture}

For $n=2$ this is the beautyful theorem of Sydler, later
reproved by Jessen and
Thorup and Dupont-Sah [DS2].

Notice that $\Omega^{0}_{\Bbb R/\Bbb Q}
= \Bbb R$, so for $i=1$ and arbitrary $n$ we get the conjecture (\ref{eucl}).

Here is a "motivic interpretation" of  conjecture (\ref{deuc}).
One could think of Euclidian geometry as of degeneration of
hyperbolic geometry. Namely,
consider the Cayley model of the hyperbolic geometry inside of the sphere
 $Q_{\varepsilon}$ given by the  equation
$\varepsilon (x_1^2 + ... x_n^2) = x_0^2$ in homogeneous coordinates.
Then in the limit
$\varepsilon \longrightarrow 0$ we will get the
Euclidean geometry in $\Bbb R^n$. Now let $\varepsilon$ be a formal
variable with $\varepsilon^2 = 0$. Then one should imagine that
$H(Q_{\varepsilon}, M):= H^n(\Bbb P^n \backslash Q_{\varepsilon}, M)$ is
 a variation (or better to say a deformation because it should split at
$\varepsilon=0$) of mixed Tate motives over $Spec \Bbb R[\varepsilon]$.
Notice that right now there is no even a hypothetical definition
of the category of mixed motives or  mixed Hodge structures over
a scheme with nilpotents. The Euclidian Dehn invariant
should have  a natural interpretation in the language of deformations of
mixed Hodge structures over $Spec \Bbb R[\varepsilon]$ similar to
the one for non euclidean Dehn invariant given in the  chapter
4.

Let
$D_T(\varepsilon)$ be the hypothetical category of mixed Tate motives over
$Spec \Bbb R[\varepsilon]$ which split over $Spec \Bbb R$.
It is natural to conjecture that
$$
Ext_{D_T(\varepsilon)}^i(\Bbb Q(0),\Bbb Q(n)) = Ker
\Bigl(gr^{\gamma}_nK_{2n-i}(\Bbb R[\varepsilon]) \longrightarrow
gr^{\gamma}_nK_{2n-i}( \Bbb R)\Bigr)\otimes \Bbb Q
$$
One can compute the right hand side using the comparision with
the cyclic homology ([FT],
[Goo]) and get the following theorem (I am indebted to
B. Tsygan for scketching  me  the   proof):
\begin{theorem}
$Ker
\Bigl(gr^{\gamma}_nK_{2n-i}(\Bbb R[\varepsilon]) \longrightarrow
gr^{\gamma}_nK_{2n-i}( \Bbb R)\Bigr)\otimes \Bbb Q =
\Omega^{i-1}_{\Bbb R/\Bbb Q}$
\end{theorem}

Therefore one should get a canonical map $H^i{{\cal
    P}}^{\bullet}_{E}(2n-1)  \longrightarrow \Omega^{i-1}_{\Bbb
  R/\Bbb Q}$. I hope  that it is an isomorphism.

 \section  {Volumes of hyperbolic $n$-manifolds and continuous cohomology
of $SO(n,1)$}

  Any oriented hyperbolic manifold $M^n$, dim$M^n = n$,
can be represented as a quotient
$M^n = \G \setminus {\cal H}^n$
where $\G$ is a torsion free subgroup of $SO(n,1)$.

\vskip 3mm \noindent
{\bf 1}.Let us denote by $BG$ the classifying space of a group $G$.
Set $G^n:= \underbrace{G\times \cdots \times G}_{n \mbox{ times}}$.
There is Milnor's simplicial model for $BG$: $BG_\bullet = EG_\bullet/G$ where
$$
EG_\bullet :\ G \leftleftarrows G^2 \
\hbox{{$\leftarrow$}\kern-10pt\lower4pt\hbox{$\leftarrow$}
\kern-13.2
pt\raise4pt\hbox{$\leftarrow$}} \
G^3 \
\hbox{{$\cdots$}\kern-12.6pt\lower4pt\hbox{$\leftarrow$}
\kern-13.5pt\raise4pt\hbox{$\leftarrow$}} \
\ldots
$$

The inclusion $j: \G \hookrightarrow SO(n,1)$ induces a map
$$
j_*:B\G \to BSO(n,1)^\de
$$
(Here $G^\de$ is a Lie group considered as a discrete group).
Notice that $\G \setminus {\cal H}^n = B\G$.

\vskip 3mm \noindent
 {\bf 2. Homological interpretation of the volume of a compact
hyperbolic manifold}.  Now let us suppose that $M^n$ is {\it compact.}
Then $H_n(M^n) = \ZZ$.  Let us denote by $b_n$ the generator of $H_n(M^n)$.
Then
 \begin{equation} \label{h2.1}
j_*(b_n) \in H_n(BSO(n,1)^\de) = H_n(SO(n,1)^\de)
\end{equation}
Further, let $G$ be a Lie group, $C^n(G)$: the space of continuous
functions on $G^n$.  There is a differential
$$
d: C^n(G) \to C^{n+1}(G) ,
df(g_0, \ldots, g_n) = \sum_{i=0}^n (-1)^i f(g_0, \ldots, \hg_i, \ldots, g_n)
$$
Let $C^n(G)^G$ is the subspace of functions invariant under the left diagonal
action of $G$ on $G^n$.  Then by definition
\begin{equation}
H_c^*(G,R) := {\cal H}^{*+1}(\ldots \stackrel {d}
{\longrightarrow}  C^n(G)^G \stackrel {d}
{\longrightarrow} C^{n+1}(G)^G
\stackrel{d}
{\longrightarrow} \ldots )
\end{equation}
are the continuous cohomology of the Lie group $G$.

Let us denote by $I(g_0x, \ldots, g_nx)$ the geodesic simplex in the hyperbolic
space ${\cal H}^n$ with vertices at points $g_0x, \ldots, g_nx$, where $g_i \in
SO(n,1)$ and $x$ is a given point in ${\cal H}^n$.  Then

\begin{equation}  \label{2.2}
\vo(I(g_0x, \ldots, g_nx))
\end{equation}
is a continuous $n$-cocycle of $SO(n,1)$ because it

a) is invariant under the left diagonal action of $SO(n,1)$

b) satisfies the cocycle condition
$$
\sum_{i=0}^{n+1} (-1)^i \vo(I(g_0x, \ldots, \widehat{g_ix},\ldots, g_{n+1}x))
= 0
$$
Let
\begin{equation}  \label{2.3}
v_n \in H_c^n(SO(n,1),R) \subset H^n(SO(n,1)^\de, R)
\end{equation}
be the cohomology class of this cocycle.

\vskip 3mm \noindent
 {\bf Theorem 3.1}. {\it Let $M^n$ be a compact hyperbolic manifold.  Then}
$$
\vo(M^n) = \langle v_n, j_*(b_n)\rangle
$$

\vskip 3mm \noindent
{\bf Proof}. There is a triangulation
$$
M^n = \bigcup_kI_k
$$
of $M^k$ in geodesic simplices $I_k$.  To obtain it choose $N$ generic points
$y_1, \ldots, y_N$ in $M^n$ and consider the corresponding Dirichlet domains

\begin{equation}  \label{2.4}
 \cd (y_j) := \{x\in M^n | \rho(x, y_\a) \le \rho(x, y_\b)
\mbox{ for any }\b \ne \a\}
\end{equation}

Then the dual triangulation is the desired one.  (The vertices of the dual
simplices are the points $y_\a$, two vertices $y_{\a_1}$ and $y_{\a_2}$
are connected by the edge if domains $\cd (y_{\a_1})$ and $\cd (y_{\a_2})$
have common codimension 1 face and so on.)

Let $\pi: {\cal H}^n \to \G\setminus {\cal H}^n$.  The group $\G$ acts on
${\cal H}^n$ freely.
Therefore $\pi^{-1}(v_\a)$ is a principal homogeneous space of $\G$.  Let
us choose elements $g_\a \in SO(n,1)$ such that $g_\a x \in \pi^{-1}(v_\a)$.
(and so $\pi^{-1}(v_\a) = \G \cdot g_\a x$).  Let $I(g_0^{(k)}x, \ldots,
g_n^{(k)}x)$ be a geodesic simplex in ${\cal H}^n$ that projects onto $I_k$.
One can choose elements $g_i^{(k)}$ in such a way that $g_i^{(k)} \in  \G \cdot
g_\a x$ for some $\a$.  Let us do this and consider the following $n$-chain
in $BSO(n,1)^\de$.

\begin{equation}  \label{2.5}
\sum_k (g_0^{(k)}, \ldots, g_n^{(k)})
\end{equation}

\vskip 3mm \noindent
 {\bf Lemma 3.2}.

a){\it The boundary of this $n$-chain is zero.

b) Its homology class in $H_n(BSO(n,1)^\de)$ coincides with}
$j_*(b_n)$.

\vskip 3mm \noindent
{\bf Proof}.
a) Let $I(h_0^{(1)}x, \ldots, h_{n-1}^{(1)}x)$ and
$I(h_0^{(2)}x, \ldots, h_{n-1}^{(2)}x)$
be codimension 1 faces of geodesic simplices
$I(g_0^{(k)}x, \ldots, g_{n}^{(k)}x)$ that project to the same face
in $\G\setminus {\cal H}^n$. Then (thanks to the special choice of elements
$g_i^{(k)}$) there is an element $\g \in \G$ such that $\g \cdot h_i^{(1)}
=  h_i^{(2)}$ for all $i$.
Further, it is clear that these faces will appear with opposite
signs.  The statement a) is proved.
The proof of statement b) follows from the definitions.

The value of the cocycle (\ref{2.2}) on the cycle (\ref{2.5}) is equal to
$\vo(\G\setminus
{\cal H}^n)$ just by  definition.  So theorem 3.1 follows immediately from
Lemma 3.2.

\vskip 3mm \noindent
 {\bf 3. The strategy for noncompact hyperbolic manifolds}.
Now let $M^n = \G\setminus {\cal H}^n$ be a noncompact hyperbolic manifold with
finite volume.  We would like to prove an analog of theorem 3.1.  The
first problem is that $H_n(M^n)=0$.
Let $H$ be a subgroup of a group $G$.  Set
$$
H_*(G,H;\QQ):=H_*(\cone(BH_\bullet\to BG_\bullet),\QQ)
$$
We will show that $M^n$ does produce an element

\begin{equation}  \label{2.6}
\tc(M_n) \in H_n(SO(n,1)^\de, T_n(s)^\de;\QQ)
\end{equation}
 where $T_n(s)$ is a subgroup of $SO(n,1)$ consisting of transformations
preserving $s\in \partial {\cal H}^n$
and acting as translations on horospheres based at $s$.
(If $s$ is the point at infinity in the upper half-space realization of ${\cal
H}^n$, then $T_n(s)$
is the group of translations in the hyperplane $x_{n}=c$ where $x_{n}$
is the vertical coordinate).

Further, let $H$ be a Lie subgroup of a Lie group $G$.  The inclusion
$i: H \hookrightarrow G$ induces a homomorphism of complexes $C^*(G)^G
\longrightarrow  C^*(H)^H$.  Set
$$
H_c^*(G,H;\RR) := H^*\cone (C^*(G)^G \stackrel {i^{*}}
{\longrightarrow}  C^*(H)^H) [-1]
$$
By definition an $n$-cocycle in $\cone (C^*(G)^G \stackrel {i^{*}}
{\longrightarrow}  C^*(H)^H)[-1]$
is a pair $(f,h)$ where $f$ is a continuous $n$-cocycle of $G$, $h$ is an
$(n-1)$-cochain on $H$ and
$$
i^*f = dh
$$
In other words cocycles in $\cone(i^*)$ are just those cocycles of $G$ whose
restriction to $H$ is cohomologous to zero.
Choose a point $x \in {\cal H}^n$.  Then $(h_i \in T_n(s))$

\begin{equation}  \label{2.7}
\tv_n(x) := (\vo I (g_0x, \ldots, g_nx), \quad \vo I(s, h_0x, \ldots,
h_{n-1}x))
\end{equation}
 is a cocycle representing certain cohomology class
$$
\tv_n \in H_c^n(SO(n,1), T_n(s))
$$
Indeed,
$$
\vo I(h_0x, \ldots, h_nx) = \sum_{i=0}^n (-1)^i \vo I (s, h_0x, \ldots,
\widehat{h_ix}, \ldots, h_nx)
$$
We will prove that $\vo(M^n) = \langle \tv_n, \tc(M^n)\rangle$.

To produce the class $\tc_n(M^n)$ we have to recall some basic facts about
noncompact hyperbolic manifolds (see \S5 of Thurston's book [Th])

\vskip 3mm \noindent
 {\bf 4. A decomposition of hyperbolic manifolds}  Let $x \in {\cal H}^n$.
Denote by $\G_\ve(x)$ the subgroup generated by all elements of $\G$ which
move $x$ to a distance $\le \ve$.

\vskip 3mm \noindent
 {\bf Theorem 3.3} (Kazhdan-Margoulis [KM]).  {\it There is an $\ve > 0$ such
that for every discrete group $\G$ of isometries of ${\cal H}^n$ and for every
$x\in
{\cal H}^n$, $\G_\ve(x)$ has an abelian subgroup of finite index.}

Let $s \in \partial {\cal H}^n$.  A part of ${\cal H}^n$ located inside an
horosphere centered
at $s$ is called horoball.
Iqn the upper half space realization of ${\cal H}^n$ horoball
centered at infinity is defined as $\{x_{n}\ge c\}$ .
Let us denote by $E_n(s)$ the stabilizer in $SO(n,1)$ of the point $s$.
Let $\De_n(s)$ be a discrete subgroup of $E_n(s)$ that has a subgroup of finite
index isomorphic to $\ZZ^{n-1}$ and acting by translation on horospheres
centered at $s$. Suppose also that $\De_n(s)$ acts freely in ${\cal H}^n$.
Then the
quotient of a horoball centered at $s$ is called cusp.

Denote by $\ell(v)$ the length of the shortest geodesic loop based at $v$.

\vskip 3mm \noindent
 {\bf Theorem 3.4}. (Decomposition Theorem). {\it  Let $M^n$ be an orientable
hyperbolic manifold of finite volume.  Then there is $\de > 0$ such that
$\ell^{-1}[0,\de] \subset M^n$ consists of finitely many components and each
of these components is isometric to a cusp.}

This theorem follows from theorem 3.3.

Recall that $M^n = \G\setminus {\cal H}^n$.  Let us denote by $C(\G)$ the
subset of
$\partial {\cal H}^n$ consisting of the points $c$ such that

a) the isotropy subgroup of  $c$ in $\G$ is nontrivial

b) the isotropy subgroup of any geodesic ending at $c$ in $\G$ is trivial

\vskip 3mm \noindent
{\bf Remark 3.5}. A nontrivial element $\g$ of $\G$ preserving a point
$c \in \partial {\cal H}^n$ and a geodesic ending at $c$ can not stabilize
points
on this geodesic because by the assumption each point of ${\cal H}^n$ has
trivial
stabilizer in $\G$.  So the element $\g$ move the geodesic along itself.
Therefore the point $c$ does not correspond to a cusp of $\G\setminus {\cal
H}^n$.
It follows from the Decomposition theorem that $C(\G)$ consists of a finite
number of $\G$-orbits. (Each of them corresponds to a cusp of
$\G\setminus {\cal H}^n$).
Let us choose $\{c_i\} \in C(\G)$ such that

\begin{equation}  \label{2.9}
C(\G) = \bigcup_i \G c_i
\end{equation}
The manifold $M^n$ is compact if and only if $C(\G) = \emptyset$.

\vskip 3mm \noindent
 {\bf 5. A triangulation of a noncompact hyperbolic manifold by ideal
geodesic simplices}.  We would like to define a Dirichlet decomposition of
${\cal H}^n$ corresponding to the set $C(\G) \in \partial {\cal H}^n$, and then
take the
dual triangulation.  The problem is that the distance $\rho(x,c)$ from a point
$x \in {\cal H}^n$ to $c \in C(\G)$ is infinite.  To define a regularized
distance
$\tro(x,c)$ let us choose a horosphere $h(c)$ for each $c \in C(\G)$ in such a
way that $h(\g \cdot c) = \g \cdot h(c)$ for any $\g \in \G$.  Then by
definition $\tro(x,c)$ is the distance from $x$ to the horosphere $h(c)$.  It
is
negative if $x$ is inside of $h(c)$.  Now the Dirichlet domains are defined as
usual:

\begin{equation}  \label{2.10}
\cd (c) :=\{x\in {\cal H}^n | \tro(x,c) \le \tro(x,c') \mbox{ for any } c' \in
C(\G),\
c'\ne c\}
\end{equation}
 They are polyhedrons with finite number of faces.

For a generic choice of horospheres $h(c_i)$ the dual polyhedrons are
simplices with vertices at points of $C(\G)$.  This triangulation is
$\G$-invariant.

\vskip 3mm \noindent
\noindent{\bf 6.}
 {\bf Theorem 3.6}. {\it Any hyperbolic manifold of finite volume $M^n$ defines
canonically a class
$$
\tc(M^n) \in H_n(SO(n,1)^\de, T_n(v)^\de;\QQ)
$$
such that}
$$
\langle \tv_n, \tc(M^n)\rangle =\vo(M^n)
$$

\vskip 3mm \noindent
{\bf Proof}.
 Let $\widetilde{\partial {\cal H}^n}$
be the set of all pairs $(v,\ell)$ where $v \in \partial {\cal H}^n$ and $\ell$
is
a geodesic ended at $v$.

Let $\widetilde{C(\G)} \subset \widetilde{\partial {\cal H}^n}$ be the set of
all pairs
$(c,\ell)$ where $c \in C(\G)$ and $\ell$ is an edge of one of geodesic
simplices from the constructed triangulation of ${\cal H}^n$. Notice that
the group $\G$ acts freely on $\widetilde{C(\G)}$
Choose geodesics $\ell_{ij}$ such that $(c_i, \ell_{ij})\subset
\widetilde{C(\G)}$ and

\begin{equation}  \label{2.11}
\widetilde{C(\G)} = \bigcup_{i,j}\G \cdot (c_i, \ell_{ij})
\end{equation}
 Take a point $(v_0, \ell_0) \in \widetilde{\partial {\cal H}^n}$ such that
$v_0 \notin
C(\G)$.  There exist elements $g_{ij}\in SO(n,1)$ with properties

\begin{equation}  \label{2.12}
a)  (c_i, \ell_{ij}) = g_{ij}\cdot (v_0, \ell_0) \qquad
b)  g_{ij}^{-1} g_{ik} \in T_n(v_0)
\end{equation}
(Recall that
$T_n(v_0)$ is the subgroup of translations).  Therefore
any $(c,\ell) \in \widetilde{C(\G)}$ can be written as

\begin{equation}  \label{2.13}
(c,\ell) = g_{(c,\ell)} \cdot (v_0, \ell_0),
\end{equation}

\begin{equation}  \label{2.14}
g_{(c,\ell)} := \g_{(c,\ell)}\cdot g_{ij},
\end{equation}
where $\g_{(c,\ell)} \in \G$ is defined from the condition $(c,\ell) =
\g_{(c,\ell)} \cdot (c_i, \ell_{ij})$

Now let us choose a representation $S_\a$ in each class of $\G$-equivalence
of simplices from the constructed triangulation of ${\cal H}^n$.  This is a
finite
set $\{S_\a\}$ and $\{\pi(S_\a)\}$ is a triangulation of $\G\setminus {\cal
H}^n$.
(Recall that $\pi: {\cal H}^n \to \G\setminus {\cal H}^n$).

Let $F(\{S_\a\}) \subset \widetilde{C(\G)}$ be the set of all pairs $(v,\ell)$
where $v$ is a vertex of a simplex $S_\a$ and $\ell$ is its edge.  Each simplex
$S_\a$ defines $n(n+1)$ elements of the set $F(\{S_\a\})$.  They correspond
to vertices of the ``truncated simplex'' $\ts_\a$.  Let
us subdivide polyhedrons $\ts_\a$ on simplices $\ts_{\a,\b}$.  Denote by
$f_{\a,\b}^0, \ldots, f^n_{\a,\b}$ the vertices of $\ts_{\a,\b}$.  According
to (\ref{2.13}) there are uniquely defined elements $g_{\a,\b}^i \in SO(n,1)$
such that

\begin{equation}  \label{2.15}
f_{\a\b}^i = g_{\a\b}^i \cdot (v_0, \ell_0) \qquad (0 \le i \le n)
\end{equation}

Set

\begin{equation}  \label{2.16}
\tc'(M^n) := \sum_{\a,\b} (g_{\a,\b}^0, \ldots, g_{\a,\b}^n) \in
BSO(n,1)(n)
\end{equation}

The boundary of chain (\ref{2.16}) is computed as follows. Each
vertex $v_a(S_\a)$ of simplex $S_\a$ gives us $n$ elements
$$
(v_a(S_\a), \ell_{a}^1(S_\a)), \ldots, (v_a(S_\a),\ell_{a}^n(S_\a))
$$
in $F(\{S_\a\})$.  So there are $h_{a}^b \in SO(n,1)$, chosen according to
(\ref{2.13})-(\ref{2.14}) such that

\begin{equation}
(v_a(S_\a),\ell_{a}^b(S_\a)) = h_{a}^b(S_\a)\cdot(v_0, \ell_0)
\end{equation}
 ($\{h_{a}^b(S_\a)\} = \{g_{\a,\b}^i\}$ of course).  So each vertex $v_a(S_\a)$
produces a chain

\begin{equation}  \label{2.18}
(h_{a}^1(S_\a), \ldots, h_{a}^n(S_\a)) \in BSO(n,1)_{(n)}
\end{equation}
It follows from the definitions that

\begin{equation}  \label{2.19}
\partial \tc_n'(M^n) = \sum_{v_a(S_\a)} (h_{a}^1(S_\a),\ldots, h_{a}^n(S_\a))
\end{equation}
where summation is over all vertices of simplices $S_\a$.

\vskip 3mm \noindent
 {\bf Lemma 3.7}. {\it For every $c_i$ (see (\ref{2.9}))

\begin{equation}  \label{2.20}
\sum_{v_a(S_\a) \in \{\G \cdot c_i\}}
(h_{a}^1(S_\a), \ldots, h_{a}^n(S_\a)) \in BSO(n,1)_{n-1}
\end{equation}
is a cycle}.

\vskip 3mm \noindent
{\bf Proof}.
Clear.

\vskip 3mm \noindent
 {\bf Lemma 3.8}. {\it The cycle (\ref{2.20}) is the image of a cycle
$\tc_{(i)}^{\prime\prime}(M^n)$ in $BG^{(i)}(v_0)_{(n-1)}$ where
$G^{(i)}(v_0)$ is a semidirect product of a finite group and the group
of translations $T(v_0)$.}

\vskip 3mm \noindent
{Proof}. Look at definitions (\ref{2.12}).  Each element $h_{a}^b(S_\a)$ is
equal to
$$
\g \cdot g_{ik} \cdot (v_0, \ell_0) = \g \cdot g_{ij} \cdot t \cdot
(v_0,\ell_0)
= g_{ij} \cdot (g_{ij}^{-1} \cdot \g \cdot g_{ij}) \cdot t \cdot (v_0, \ell_0)
$$
where $\g \in \De_n(c_i)$ (the stabilizer of $c_i$ in $\G$) and $t \in T(v_0)$.
Notice that $g_{ij}^{-1}T(c_i) g_{ij} = T(v_0)$ and $\De_n(c_i)$ is a
semidirect
product of a finite group and a subgroup in $T(c_i)$.  Lemma 3.8 is proved.

It follows from Lemmas 3.7 and 3.8 that

\begin{equation}  \label{2.21}
c_E(M^n) := \left( \sum_i \tc_{(i)}^{\prime\prime}(M^n), \tc'(M^n)\right)
\end{equation}
 is a cycle in
$$
\cone(BG^{(i)}_n(v_0)_\bullet \to BSO(n,1)_\bullet)
$$

\vskip 3mm \noindent
 {\bf Proposition 3.9}. $\langle \tv_n, c_E(M^n) \rangle = \vo(M^n)$.

\vskip 3mm \noindent
{\bf Proof}.  The cocycle
\begin{equation}  \label{2.22}
v_n(x) = \left( \vo I(g_ox, \ldots, g_nx), \vo I(v_o, h_0x, \ldots, h_{n-1}x)
\right)
\end{equation}
represents the cohomology class $\tv_n \in H_c^n(SO(n,1), E_n(v_0))$.
The cocycles $v_n(x)$ for different points $x \in {\cal H}^n$ are (canonically)
cohomologous.  Now let us move point $x$ in (\ref{2.22}) to the boundary point
$v_0$.
Then the limit of value of the component $\vo (I(g_0x, \ldots, g_nx))$ on the
chain $\tc'(M^n)$ (see (\ref{2.16})) exists and is equal to $\sum_\a \vo(S_\a)
= \vo
(M^n)$.  From the other hand the limit of the second component $\vo I(v_0,
h_0x, \ldots, h_{n-1}x)$  of cocycle $v_n(x)$ on the chain $\sum_i c_{(i)}^{
\prime \prime}(M^n)$ is zero.  Proposition 3.9 is proved.

To prove theorem 3.6 it remains to show that there is a cycle $c_T(M^n)$ in
$$
\cone(BT_n(v_0)_\bullet \to BSO(n,1)_\bullet)
$$
homologous to $N \cdot c_E(M^n)$ for certain integer $N$.
This is a consequence of the following

\vskip 3mm \noindent
 {\bf Lemma 3.10}.  {\it The homomorphism
$$
j_*: H_*(T_n(v_0),\QQ) \to H_*(G^i(v_0), \QQ)
$$
induced by the inclusion $j: T_n(v_0) \hookrightarrow G^{(i)}(v_0)$
is a map onto}.

\vskip 3mm \noindent
{\bf Proof}.  Set $A:= G^{(i)}(v_0) / T_n(v_0)$.  There is
Hochshild-Serre spectral sequence
$$
E_{pq}^2 = H_p(A, H_q(T_n(v_0), \QQ)) \Rightarrow H_{p+q}(G^{(i)}(v_0),\QQ)
$$
Further, for a finite group $A$ and an $A$-module $V$,
$H_i(A, V\otimes \QQ) = 0$ for $i>0$.  Therefore
$$
H_*(G^{(i)}(v_0), \QQ) = H_*(T_n(v_0),\QQ)_A.
$$
Lemma 3.10 and hence theorem 3.6 is proved.

\vskip 3mm \noindent
 {\bf 7. $\tc(M^{2n-1})$ and $H_{2n-1}(GL_{N}(\bar {\Bbb Q}),\Bbb Q)$}.
It is usefull to keep in mind the following trivial lemma

{\bf Lemma 3.11} {\it Let $\tilde G \rightarrow G$ be a surjective
homomorphism of groups with finite kernel. Then it induces isomorphism on
homology with rational coefficients.}

Let $N=2^{n-1}$ and $s_+: Spin(2n-1,1) \rightarrow GL_{N}(\Bbb
C)$ be the halfspinor
representation. Set $\tilde T_{2n-1}(v):= s_+ (T_{2n-1}(v))$. We
will produce a
basis in this representation and
a 1-dimensional subgroup $H \subset
diag(GL_N)$ isomorphic to  ${\Bbb G}_m$ and diagonal in this basis such
that the following 3 conditions hold.

1. $s_+(Lie(T_{2n-1}(v)))$ is an abelian Lie algebra contained in the upper
triangular Lie algebra.

2. It is normalized by a subgroup $H$.

3. The eigenvalues of an element $t \in H(\Bbb Q) = {\Bbb Q}^*$ acting on
$s_+(Lie(T_{2n-1}(v)))$ by conjugation are equal to $t$.

For this  let me remind the construction of the
halfspinor representation.
The spinor representation of the complex Lie algebra $o(V_{2n})$, where
$V_{2n}$
is a $2n$-dimensional complex vector space,  can
be described as follows. Let us choose a decomposition $V_{2n} = U
\oplus W$ where $U$ and $W$ are (maximal) isotropic planes.
Then $o(V_{2n})$ is isomorphic to a Lie subalgebra of the Lie algebra of all
superdifferential
operators of order $\leq 2$ acting on the vector space $\Lambda^{\ast}W$
(which is the space of regular functions on the odd variety $W$). To be
precise, $o(V_{2n})$ is the Lie subalgebra of all superdifferential
operators of total  degree $0$ or $2$.

Notice that
$$
\Lambda^{\ast}W = S^- \oplus S^+
$$
where $S^- = \Lambda^{odd}W$, $S^+ = \Lambda^{even}W$.  Each subspace is
preserved by the action of $o(V_{2n})$. The
corresponding representations are the half spinor representations $s_+$ and
$s_-$.

Let me describe the spinor representation in
coordinates. Choose a basis $e_1,...,e_n;f_1,...,f_n$ in $V_{2n}$ such that
$W=<e_1,...,e_n>; U=<f_1,...,f_n>$ and $(e_i,f_j)=\delta_{ij}$. Then an
element $X \in o(2n)$ looks in this basis as follows:
$$
\begin{array}{cc}
 A & B \\
 C & -A^t
\end{array}
$$
where $B^t =-B$ and $C^t =-C$. One has
$$
s: X \longmapsto  - \frac{1}{2} tr A + \sum_{1 \leq i,j \leq n}
a_{ij}\xi_i \partial_{\xi_j} +
\sum_{1<i<j<n}\Bigl( b_{ij}  \xi_i \wedge
\xi_j  + c_{ij}\partial_{\xi_i} \wedge \partial_{\xi_j}\Bigr)
$$

Let ${\cal N} \subset s(o(2n))$ be the Lie subalgebra
 acting
in the halfspinor representation  by the following operators:
$$
\sum_{1<j<n}a_{1j}\xi_1 \partial_{\xi_j} + \sum_{1<j<n}b_{1j} \xi_1
\wedge \xi_j
$$
We will denote by  $N$ the  corresponding subgroup  of the group $GL_N$.
The group
 $\tilde T_{2n-1}(v)$
 is conjugated to a subgroup of  N.

Let us embed the group $\Bbb G_m$ to
$GL_N$ as follows:   $ t \in \Bbb G_m$ multiplies a
 vector $\xi_{i_1} \wedge ... \wedge \xi_{i_k}$ where $i_1 <
 ... < i_k$ by $t$  if
 $i_1 =1$ and by $1$ otherwise. This is the subgroup  $H$.
Let
$P$ be the semidirect product of the groups $H$ and $N$.

Let us denote by $\hc(M^{2n-1})$ the image of $\tc(M^{2n-1})$ under the
composition of natural
homomorphisms
$$
H_{2n-1}(SO(2n-1,1),  T_{2n-1}(v); \Bbb Q) \stackrel {{s_+}_{\ast}}
{\longrightarrow}  H_{2n-1}(GL_{N}(\Bbb C), \tilde T_{2n-1}(v);\Bbb Q)
\longrightarrow
$$
$$
H_{2n-1}(GL_{N}(\Bbb C), P(\bar {\Bbb Q});\Bbb Q)
$$
The second arrow is provided by a conjugation in $GL_N(\bar {\Bbb
  Q})$.
Let $P(\bar {\Bbb Q})$ be the group of $\bar {\Bbb Q}$-points of $P$.

\vskip 3mm \noindent
{\bf Lemma 3.12}
a) {\it Let $M^{2n-1}$ be a compact hyperbolic manifold. Then
  (see \ref{h2.1})
$$
j_{*}(b_{2n-1}) \in H_{n}(SO(2n-1,1)(\bar
  {\Bbb Q}), \Bbb Z)
$$

b)Let $M^{2n-1}$ be a noncompact hyperbolic manifold. Then}
$$
\hc(M^{2n-1}) \in H_{2n-1}(GL_{N}(\bar {\Bbb Q}), P(\bar {\Bbb Q})
 ;\Bbb Q)
$$

\vskip 3mm \noindent
{\bf Proof} According to  Weil's rigidity, see 6.6, 6.7, 7.13 in [Ra],
the subgroup $\G \in  SO(2n-1,1)$ is conjugate to a subgroup whose entries
are algebraic numbers.

There is the usual exact sequence of cone:
$$
\ldots \to
H_{2n-1}(P(\bar {\Bbb Q});\Bbb Q) \longrightarrow
H_{2n-1}(GL_{N}(\bar {\Bbb Q});\Bbb Q) \stackrel {e_{\ast}}
{\longrightarrow}
H_{2n-1}(GL_{N}(\bar {\Bbb Q}),P(\bar {\Bbb Q});\Bbb Q)
$$
\begin{equation}
\longrightarrow
H_{2n-2}(P(\bar {\Bbb Q});\Bbb Q) \longrightarrow
H_{2n-2}(GL_{N}(\bar {\Bbb Q});\Bbb Q)
\end{equation}

\vskip 3mm \noindent
 {\bf Theorem 3.13}.  {\it There exists an element
$c(M^{2n-1}) \in H_{2n-1}(GL_{N}(\bar {\Bbb Q});\Bbb Q)$ such
that} $\hc(M^{2n-1}) = e_* c(M^{2n-1})$.

\vskip 3mm \noindent
{\bf Proof}.
It follows from

 {\bf Proposition 3.14}. {\it $i_*: H_*(P(\bar {\Bbb Q});\Bbb Q) \to
H_*(GL_{N}(\bar {\Bbb Q}); \Bbb Q)$ is injective.}

{\bf Proof}. Recall that  $P = H\cdot N$ where $H$ is the
diagonal part of $P$ and
$N$ is the unipotent one. Notice that $H(\bar {\Bbb Q}) = \bar {\Bbb Q}^*$.

\vskip 3mm \noindent
  {\bf Lemma 3.15}. {\it The natural map
$H_*(H(\bar {\Bbb Q});\Bbb Q) \to H_*(P(\bar {\Bbb Q});\Bbb
Q)$ is an isomorphism.}

\vskip 3mm \noindent
{Proof}.  There is the Serre-Hochshild spectral sequence
$$
H_p(H,H_q(N,\Bbb Q)) \Rightarrow H_{p+q}(P;\Bbb Q).
$$
For a number field $F$ the group $N(F)$
is isomorphic to the additive group
of a finite dimensional $\Bbb Q$-vector
space $V(F)$.  So $H_q(N(F),\Bbb Q) \simeq \land_{\Bbb Q}^q V(F)$.
Further, any integer $n
\in \Bbb Q^*
=  H(\Bbb Q)$ acts on $H_q(N,\Bbb Q)$ by multiplication on $n^q$.
Therefore $H_*(H(\bar {\Bbb Q}),H_q(N(\bar {\Bbb Q});\Bbb Q))$
is anihilated by multiplication on $n^q-1$.

This is a particular case of the following general fact (Proposition 5.4 in
[M]): let $G$
be a group and $V$ is a $G$-module.  Then the action of any element $g\in G$ on
$G$ by
conjugation and on $V$ via the $G$-module structure induces the identity map on
$H_*(G,V)$.
Lemma 3.15 is proved.

To complete the proof of the proposition 3.14 consider the homomorphism
$det: GL_N({\bar {\Bbb Q}}) \longrightarrow \bar \Bbb
Q^{\ast}$. Notice that
${\bar {\Bbb Q}}^{\ast} = H({\bar {\Bbb Q}}) \rightarrow GL_N({\bar {\Bbb Q}})
\stackrel{det}{\rightarrow} {\bar {\Bbb
Q}}^{\ast}$ is given by $x \longmapsto x^k$ for a certain positive integer $k$.
Therefore its composition with  $H_{\ast}(H({\bar {\Bbb
Q}});\Bbb Q) \longrightarrow H_{\ast}(GL_N({\bar {\Bbb Q}});\Bbb Q)
\longrightarrow H_{\ast}({\bar {\Bbb Q}}^{\ast};\Bbb Q)$ is injective.
The proposition is proved.

\vskip 3mm \noindent
{\bf 8. The spinor representation and the Pfaffian}.
The supertrace $Str$ in the spinor representation $S = S^+ \oplus
S^-$ is defined as follows:
$$
Str:= Tr|S_+ - Tr|S_- \qquad \mbox{i.e.} \quad Str(A):= Tr s_+(A)
- Tr s_-(A)
$$

There is a remarkable invariant polynomial of degree $n$ on the Lie
algebra $o(2n)$: the Pfaffian.
Its
restriction to the Cartan subalgebra $(t_1,...,t_n)$,
$e_i \longmapsto t_i e_i , f_i \longmapsto - t_i f_i$
is given by formula $Pf(t_1,...,t_n) = t_1 \cdot ... \cdot t_n$.
For  a skewsymmetric $2n \times 2n$ matrix $A$  one has $Pf(A)^2=
(-1)^n det (A)$.

{\bf Proposition 3.16} {\it The invariant polynomial $A
  \longmapsto (-1)^{n(n-1)/2}n! \times Pf(A)$  on $o(2n)$
is equal to the  restriction of
the invariant polynomial
$A \longmapsto Str(A^n)$ in  the spinor representation}.

{\bf Proof}. Let us compute the restriction of the invariant polynomial
$Str(A^n)$ in the spinor representation to the Cartan subalgebra
$
e_i \longmapsto t_i e_i , f_i \longmapsto - t_i f_i
$

We will get
$$
Str\Bigl( \sum_{i=1}^n t_i \xi_i \partial_{\xi_i} -
\frac{1}{2}(t_1+...+t_n)\Bigr)^n = (-1)^{n(n-1)/2}n!\cdot t_1\cdot
... \cdot t_n
$$
Indeed, it is easy to see that the supertrace of the superdifferential
operator $\xi_{i_1}\wedge
... \wedge \xi_{i_m} \partial_{\xi_1} \wedge ... \wedge
\partial_{\xi_{i_m}}$ is $0$ if
$m<n$ and
$(-1)^{n(n-1)/2}n!$ if $m=n$.

Consider the involution $\sigma$ of the Lie algebra $o(2n,C)$ which corresponds
to the involution interchanging the  ``horns'' of the Dynkin diagram,
whose involutive subalgebra is $o(2n-1,C)$.

{\bf Lemma 3.17} {\it The involution $\sigma$ of the Lie algebra $o(2n,C)$
interchanges the halfspinor
representations}.

{\bf Proof}. The
halfspinor representation is the fundamental representation corresponding
to a vertex of a ``horn''.

\vskip 3mm \noindent
{\bf 9.  The Lobachevsky class $v_{2n-1} \in
H_c^{2n-1}(SO(2n-1,1), \Bbb R)$ can be obtained by restriction of
the Borel class in the halfspinor representation}.
Recall that for a  Lie  group $G$  with maximal
compact subgroup $K$ one has the Van Est isomorphism
$$
H^{\ast}_c(G,\Bbb R) = H^{\ast}(Lie G, Lie K)
$$
(On the right we have the relative Lie algebra cohomology).
{}From the other hand
$$
H^{\ast}(Lie G, Lie K) = H^{\ast}(Lie G \otimes_{\Bbb
R} \Bbb C, Lie K \otimes_{\Bbb
R} \Bbb C)
$$
Let us suppose that $K = G^{\sigma}$ where $\sigma$ is an involution
 of the group $G$, (i.e. $G,K,\sigma)$ is an involutive pair).
Let $G_{\Bbb C}$ (resp $K_{\Bbb C}$) be the complex Lie group corresponding
to $Lie G \otimes_{\Bbb R} \Bbb C$ (resp $Lie K \otimes_{\Bbb R} \Bbb C$).
Then $K_{\Bbb C}$ is the fixed point set of an involution $\sigma_{\Bbb C}$
of the complex Lie
group $G_{\Bbb C}$.
Finally, let $G_u$ be the maximal compact subgroup of $G_{\Bbb C}$.
 Then there is an
involution $\sigma_u$ of $G_u$ such that $K = G_u^{\sigma_u}$.
In this situation on has the isomorphism
$$
H^{\ast}_c(G,\Bbb R) = H^{\ast}(Lie G, Lie K) = H_{top}^{\ast}(G_u/K, \Bbb R)
$$
which is functorial with respect to the maps of the symmetric pairs.

For example when $G=Spin(2n-1,1)$ one has $K= Spin(2n-1)$, $G_u = Spin(2n)$
and
$$
H^{\ast}_c(Spin(2n-1,1),\Bbb R) = H^{\ast}_c(SO(2n-1,1),\Bbb R) =
H_{top}^{\ast}(S^{2n-1}, \Bbb R)
$$
So there is just one (up to a scalar) nontrivial continuous cohomology class of
the Lie
group $SO(2n-1,1)$. Such a class $v_{2n-1}$ was produced in s.2 above.
We will call it the Lobachevsky class.

Similarly
$$
H_c^{*}(GL_{N}(\Bbb C),\Bbb R) = H_{top}^{\ast}(U(N) \times
U(N)/U(N)_{diag}, \Bbb R) =
\Lambda ^{*}(b_{1},b_{3}, ... ,b_{2N+1})
$$
where $b_{2i-1} \in H_c^{2i-1} (GL_{N}(\Bbb C),\Bbb R)$ are the
continuous cohomology classes corresponding to the primitive
generators $B_{2i-1}$ of
$H_{top}^{\ast}(U(N), \Bbb Z)$. Namely  there is
canonical projection  $p:
U(N) \longrightarrow S^{2N-1} \subset \Bbb C^N$ and  $B_{2N-1}:=
p^*[S^{2N-1}]$. Further,
the restriction map
$H^{2i-1}(U(N),\Bbb Z) \longrightarrow H^{2i-1}(U(i),\Bbb Z)$
sends $B_{2i-1}$ to the described above generator. I will
call them the Borel classes. (Borel [Bo2] used a different
normalisation of these classes).

\vskip 3mm \noindent
{\bf Theorem 3.18}. {\it Let $s_+: Spin(2n-1,1) \rightarrow Aut(S_+)$ be the
halfspinor representation. Then the restriction of the Borel class
$b_{2n-1}$ is proportional to the Lobachevsky class:
$s_+^{\ast}b_{2n-1} = c\cdot v_{2n-1}$, $c \not = 0$}.

{\bf Proof}. We will need the following 3 symmetric pairs (Lie group, an
involutive Lie subgroup)
$$
(Spin(2n-1,1), Spin(2n-1)), \quad (Spin(2n,\Bbb C),
Spin(2n-1,\Bbb C)),
$$
$$
\quad (Spin(2n),
Spin(2n-1))
$$
We start with the half spinor representation providing us with a morphism of
symmetric  pairs
$$
(Spin(2n-1,1),Spin(2n-1)) \longrightarrow (GL_{2^{n-1}}(\Bbb C), U(2^{n-1}))
$$
Then we complexify it, getting a morphism of symmetric pairs
$$
(Spin(2n, \Bbb C),Spin(2n-1, \Bbb C)) \longrightarrow
(GL_{2^{n-1}}(\Bbb C) \times GL_{2^{n-1}}(\Bbb C), GL_{2^{n-1}}(\Bbb C))
$$
Finally  restricting it to the maximal compact subgroups we get a
morphism of symmetric pairs
$$
a: (Spin(2n),Spin(2n-1)) \longrightarrow (U(2^{n-1}) \times U(2^{n-1}) ,
U(2^{n-1})_{diag})
$$
where $U(2^{n-1})_{diag}$ is the diagnal subgroup. The last map is
 $a: u \longmapsto (u,\sigma(u))$ where $\sigma$ is the
involution whose involutive subgroup is $Spin(2n-1)$.

The projections of $a (Spin(2n))$ to the first and
second factors are halfspinor and antyhalfspinor representations.
Indeed, the involution $\sigma$ of the Lie algebra $o(2n,C)$
whose involutive subalgebra is $o(2n-1,C)$ interchanges the halfspinor
representations (lemma 3.17).


So we get a commutative diagram
\begin{equation} \label{d}
\begin{array}{ccc}
Spin(2n)& \stackrel{a}{\longrightarrow} & U(N) \times U(N)\\
&&\\
\downarrow p_1& & \downarrow p_2\\
&&\\
\frac{Spin(2n)}{Spin(2n-1)} & \stackrel{a'}{\longrightarrow}
& \frac{U(N) \times
U(N)}{U(N)_{diag}}
\end{array}
\end{equation}

Under the map $p_2: (u_1,u_2) \rightarrow u_1u_2^{-1}$ the primitive
generator $B_{2n-1} \in H^{2n-1}_{top}(U(N))$ goes to $p_2^{\ast}B_{2n-1} =
B_{2n-1} \otimes 1 - 1 \otimes B_{2n-1}$

The key topological statement is  that $a'(B_{2n}) \not = 0 $.
It is a corollary of proposition 3.16. To see this
let us remind that there is  the Euler class
$$
E_{2n-1} := \pi^*([S^{2n-1}]) \in H_{top}^{\ast}(SO(2n))
$$
 where $\pi:
SO(2n) \rightarrow S^{2n-1}$ is the  natural
projection.
There are the corresponding classes
$$
e_{2n-1} \in H^{\ast}(o(2n, \Bbb C)) \qquad \mbox{and} \quad {\bar e}_{2n} \in
H_{top}^{\ast}(BSO(2n))
$$
One can identify $H_{top}^{\ast}(BSO(2n))$
with the ring of invariant polynomials $H^{\ast}(o(2n))^{o(2n)}$ on the Lie
algebra $o(2n)$. Then the Euler class corresponds to
(a non zero multiple of)
the  Pfaffian. The constant $c_n$ can be calculated as follows.
Set $E_{2n-1} = \alpha_n [Pf], B_{2n-1} = \beta_n [Tr A^n]$ where
$[Pf]$ and $[Tr A^n]$ are the topological cohomology classes of
$S^{2n-1}$ and $U(N)$ coresponding to the invariant symmetric
polynomials given by the Pfaffian and $TrA^n$. Then $c_n =
(-1)^{n(n-1)/2}n! \cdot \beta_n/\alpha_n$.

Theorem 3.18 is proved.

Theorem 1.1 follows immediately from theorems 3.18, 3.1 (in the compact
case), 3.6 (in the noncompact case) and 3.13 and lemma
3.12.
First of all the image of the element $c(M^{2n-1}) \in
H_{2n-1}(GL_N(\bar {\Bbb Q});\Bbb Q)$ in $K_{2n-1}(\bar {\Bbb Q})$ gives the
invariant of $M^{2n-1}$ we are looking for. It is well defined
because  $H_{2n-1}(P(\bar {\Bbb Q});\Bbb Q) = \bar {\Bbb
  Q}^*$ projects to
zero in $K_{2n-1}(\bar {\Bbb Q})$ (notice that $2n-1
>1$). Further, we have an embedding
$$
\varphi: \Bigl(SO(2n-1,1)(\bar {\Bbb Q}), T_n(v)(\bar {\Bbb
  Q})\Bigr) \hookrightarrow \Bigl( GL_N(\bar
{\Bbb Q}), P(\bar {\Bbb Q})\Bigr)
$$
The restriction of the cohomology class $b_{2n-1}$ to $P(\bar {\Bbb Q})$ is
zero, so it provides a relative class ${\tilde b}_{2n-1}$. Similarly one
has the relative Lobachevsky class ${\tilde v}_{2n-1}$.
It follows from theorem 3.18 that  $\varphi^*({\tilde b}_{2n-1}) = {\tilde
v}_{2n-1}$. According to theorems 3.1,  3.18 and 3.13
$$
vol(M^{2n-1}) \quad = \quad <{\tilde v}_{2n-1}, {\tilde
  c}(M^{2n-1}> \quad = \quad <{\tilde b}_{2n-1}, \varphi_*{\tilde
  c}(M^{2n-1})> =
$$
$$
\quad <b_{2n-1},
\varphi_*c(M^{2n-1})>
$$
Theorem 1.3 now follows from theorem 1.1 and main results of [G1] or
[G2] where the Borel regulator on $K_5(F)$ where $F$ is a number field was
computed via the trilogarithm. See  s. 3 and 5 in [G2].

Theorem 1.2 follows, for example, from theorem 1.1 and results of s.2 in [G1].

{\bf Remark}. The halfspinor representations
seems to be the only ones among the fundamental representations
of $SO(2n-1,1)$ such that the restriction of the Borel class to
 $SO(2n-1,1)$ is not zero. This means that only the halfspinor
 representations lead to a nontrivial invariant  $c(M^{2n-1}) \in
H_{2n-1}(GL_N(\bar {\Bbb Q});\Bbb Q)$.

{\bf 10. A plan of the second  approach to the main results}
Let us say that a geodesic simplex is defined over  $\bar
\Bbb Q$ if it is equivalent under the motion group to a
geodesic simplex which has  vertices with coordinates   in $\bar
\Bbb Q$.

{\bf Proposition 3.19} {\it Any hyperbolic manifold can be decomposed on
geodesic simplices defined over  ${\bar {\Bbb Q}}$}.

{\bf Corollary 3.20} {\it The scissor congruence class of a hyperbolic
manifold belongs to ${\cal P}({\cal H};{{\bar {\Bbb Q}}})$}.

{\bf Proof}. By Weil's rigidity theorem ([Ra]) $\pi_1(M)$ is
congugated to a subgroup whose entries are algebraic numbers.
If $M$ is compact let us take sufficiently many points whose coordinates
are algebraic
numbers inside of the fundamental domain. Consider the Dirichlet
decomposition corresponding to these points and their orbits under the
action of the group  $\pi_1(M)$. The intersection of
$n+1$ geodesic hyperplane defined over ${\bar {\Bbb Q}}$ is a defined over
${\bar {\Bbb Q}}$ point. So the vertices of this decomposition defined over
${\bar {\Bbb Q}}$.

If $M$ has cusps then their coordinates are also defined over $\bar \Bbb
Q$, and the Dirichlet decomposition with respect to chosen points and
the cusps gives the desired decomposition of $M$ on geodesic simplices
defined over ${\bar {\Bbb Q}}$ (see also [EP]).

After this applying results of the next chapter and chapter 5 we get
another, more conceptual proof of the main theorem.

\section {Noneuclidean  polyhedrons
 and mixed Tate motives}

\vskip 3mm \noindent
 {\bf 1. The mixed Hodge structure corresponding to a non-Euclidian
geodesic simplex}.
Let $Q \subset \CC P^{2n-1}$ be a nondegenerate quadric and $M=(M_1,
\ldots, M_{2n})$ be a simplex in generic position with respect to $Q$.
Recall that $h(Q,M)$ is a notation for the mixed Hodge-Tate structure  in
$H^{2n-1}({\Bbb C}{\Bbb P}^{2n-1} \setminus Q, M; \Bbb Q)$ and the weights are
from $0$ to $2n$.
There are canonical isomorphisms
\begin{equation}
\Bbb Q(0) = H^{2n-1}({\Bbb C}{\Bbb P}^{2n-1}, M; \Bbb Q) \longrightarrow
W_{0}H^{2n-1}
( {\Bbb C}{\Bbb P}^{2n-1}\setminus Q,M; \Bbb Q)
\end{equation}
\begin{equation}
gr^{W}_{2n}H^{2n-1}({\Bbb C}{\Bbb P}^{2n-1} \backslash Q, M; \Bbb Q)
\longrightarrow H^{2n-1}({\Bbb C}{\Bbb P}^{2n-1}\setminus Q; \Bbb Q) = \Bbb
Q(n)
\end{equation}

The cohomology class $[w_Q]$  of (defined up to a sign) meromorphic
$2(n-1)$ - form
\begin{equation}  \label{4.1}
w_Q := 1 /{(2\pi i)^n} \cdot \sqrt{\det \tq} \cdot\frac{\sum_{i=1}^{2n} (-1)^i
x_i dx_i
\land \ldots \land \wh{dx_i} \land \ldots \land dx_{2n}}{\tq^n}
\end{equation}
has the following geometrical
interpretation.  The quadric $Q$ has 2 families of maximal isotopic
subspaces.  Let us denote by $\a_Q$ and $\b_Q$ the corresponding cohomology
classes
in $H^{2n-2}(Q)$.

The exact sequence
$$
\ldots \to H^{2n-1}({\Bbb C}{\Bbb P}^{2n-1}\setminus Q) \longrightarrow
H_Q^{2n}
({\Bbb C}{\Bbb P}^n)
= H^{2n}(Q)(-1) \longrightarrow  H^{2n}({\Bbb C}{\Bbb P}^n) \to \ldots
$$
shows that $\a_Q - \beta_Q \in \im \de$.  In fact $\a_Q - \beta_Q = \pm
\de([w_Q])$.
So to choose a sign in (\ref{4.1}) one has to choose  one of the families of
maximal isotropic
subspaces in $Q$.

Further, an orientation of the simplex $M$ (provided by the numeration of
hyperplanes
$M_i$) corresponds to a generator of $H_{2n-1}({\Bbb C}{\Bbb P}^{2n-1}, M)$.
 Therefore if we choose a sign of ${w_Q}$ and a simplex $M$ then the mixed
Hodge structure $h(Q,M)$ has the following $n$-frame: in $gr^W_{2n}h(Q,M)$ a
vector is distinguished: $[w_Q]$, and in $gr^W_0 h(Q,M)$ a functional is
distinguished: a generator of
$H_{2n-1}(\CC P^{2n-1},M)$.

\vskip 3mm \noindent
 {\bf 3. Hodge-Tate structures and the Hopf algebra $\cal H_{\bullet}$}.
For the convenience of the reader I will just reproduce some definitions
from [BGSV].  Let us call a Hodge-Tate structure a mixed $ \Bbb Q$
Hodge structure with weight factors isomorphic to $\Bbb Q(k)$.  Say
that $H$ is $n$-framed Hodge-Tate structure if it is supplied with a
nonzero vector
in $gr^W_{2n}H$ and a nonzero functional in $gr^W_0 H$.  Consider the coarsest
equivalence
relation on the set of all $n$-framed Hodge-Tate structures for which $H_1 \sim
H_2$
if there is a map $H_1 \to H_2$ compatible with frames.  For example any
$n$-framed
Hodge-Tate structure is equivalent to a one $H$ with $W_{-2}H = 0$,
$W_{2n}H=H$.
Let ${\cal H}_{n}$ be the set of equivalence classes.

One may introduce on ${\cal H}_{n}$ the structure of an abelian group as
follows.  For
$[H_1], [H_2] \in {\cal H}_{n}$ a vector in $gr^W_{2n}(H_1 \oplus H_2)$ will be
the one
whose components are distinguished vectors in $gr_{2n}^W H_1$ and
$gr_{2n}^WH_2$.

The frames in $gr_{0}^W H_i$ define the maps $gr_{0}^W (H_1\oplus H_2) \to
\ZZ(0) \oplus \ZZ(0)$.  Its composition with the sum map $\ZZ(0) \oplus \ZZ(0)
\to \ZZ(0)$ will be a distinguished functional in $gr_{0}^W H$.

Let $-[H]$ be the class of $H$ in which the frame in $gr_{2n}^WH$ is multiplied
by $-1$.  (We will get the same class multiplying by $-1$ the functional in
$gr_0^WH$).

The tensor product of mixed Hodge structures induces the commutative
multiplication
$$
\mu: {\cal H}_{k} \otimes {\cal H}_{\ell}\to {\cal H}_{k+ \ell}
$$
Let us define the comultiplication
$$
\nu = \bigoplus_{k+\ell=n} \nu_{k \ell}: {\cal H}_{n}\to \bigoplus_{k + \ell =
n}
{\cal H}_{k}\otimes
{\cal H}_{\ell}
$$

Let $[H] \subset {\cal H}_{n}$, and $(v_n, f_0)$ is a framing for $H:v_n \in
gr_{2n}^WH$,
$f_0 \in (gr_{0}^W H)^*$.  Let $R \subset gr_{2k}^W H$ be the lattice and $R^*
\subset (gr_{2k}^W H)^*$ is the dual one.  Define homomorphisms
$$
\vp: R \to {\cal H}_{k},\qquad \psi: R^* \to {\cal H}_{n-k}
$$

Namely, for $x \in R$ (resp. $y \in R^*$), $\vp(x)$ (resp. $\psi(y)$)
is the class of the mixed Hodge structure $H$ with framing $(x,f_0)$
(resp. $(v_n,y)$).

Let $\{e_j\}$, $\{e^j\}$ be dual bases in $R$, $R^*$.  Then
$$
\nu_{k,n-k}([H]) := \sum_j \vp(e_j)\otimes \psi(e^j)
$$
It is easy to see that
$$
\nu(a\cdot b) =\nu(a)\cdot\nu(b)
$$
Where $a\cdot b = \mu(a,b)$.  Therefore the abelian group ${\cal
H}_{\bullet}:=\oplus {\cal H}_n$, where
${\cal H}_0 := \ZZ$, has a structure of graded Hopf algebra with the
commutative multiplication
$\mu$ and the comultiplication $\nu$.

The following result is the main reason to consider the Hopf algebra
${\cal H}_{\bullet}$. More detailes can be found in [BMS].

\vskip 3mm \noindent
{\bf Theorem 4.1}. {\it The category of mixed $\Bbb Q$-Hodge-Tate structures is
canonically equivalent to the category of finite-dimensional graded ${\cal
H}_{\bullet}$-comodules.

Namely, the equivalence assigns to a Hodge structure $H$ the graded comodule
$M(H)$, $ M(H)_{n} = gr^{W}_{2n}(H)$ with ${\cal H}_{\bullet}$- action $M(H)
\otimes M(H)^{\ast} \longrightarrow {\cal H}_{\bullet}$ given by the formula
$x_{m} \otimes y_{n} \longrightarrow $ class of mixed Hodge structure $H$
framed by $x_{m}, y_{n}$. }

\vskip 3mm \noindent
{\bf 3. The Hopf algebra $S(\Bbb C)_{\bullet}$. }
Let $(Q, \alpha)$ be a nondegenerate quadric $Q \in {\Bbb C}{\Bbb P}^{2n-1}$
together
with a choice of one of the families of maximal isotropic subspaces on $Q$
denoted by $\alpha$.

Let us denote by $S_{2n}(\Bbb C)$ the abelian group generated by generic pairs
$[(Q,\alpha);M]$ subject to the following relations:

R0) If $M_{0} \cap \ldots \cap M_{2n-1} \neq \emptyset$ then
$[(Q,\alpha);M] = 0$

R1) ({\it Projective invariance}). For any projective transformation $g \in
PGL_{2n}(\Bbb C)$
$$
[(Q,\alpha);M] = [(gQ,g\alpha);gM]
$$

R2) ({\it Skew symmetry}) a) For any permutation $\sigma$ of the set \{0, ...
,2n-1\}   one has
$$
[(Q,\alpha);M] = (-1)^{| \sigma|} [(Q,\alpha);M_{\sigma}]
$$
where $ M_{ \sigma} = (M_{\sigma(0)}, ... ,(M_{\sigma(2n-1)})$

b) Let $\beta$ be another family of maximal isotropic subspaces on the quadric
$Q \in {\Bbb C}{\Bbb P}^{2n-1}$   . Then
$$
[(Q,\alpha);M] = -[(Q,\beta);M]
$$

R3) ({\it Additivity}). Let $M_{0},...,M_{2n}$  be hyperplanes such that $Q
\cap M_{I}$  is a nondegenerate quadric for any subset $I \in \{0, ... ,2n\}$ .
Set
$M^{(j)} : = (M_{0},..., \hat M_{j}, ... ,M_{2n})$. Then
$$
 \sum_{j = 0}^{2n+1}(-1)^{j}[(Q,\alpha);M^{(j)}] = 0
$$


For example, $S_{2}(\Bbb C)$  is generated by 4-tuples of points
$(q_{0},q_{1},m_{0},m_{1})$ on ${\Bbb C}{\Bbb P}^{1}$. Here by definition
$\alpha  = q_{0}$. The cross-ratio defines an isomorphism $r: S_{2}(\Bbb C)
\rightarrow {\Bbb C}^{\ast}$.

\vskip 3mm \noindent
{\bf Remark 4.2}. The ordering of points and relation (R2a) is necessary
to have this isomorphism. Moreover if we have it, the orientation of the
quadric $Q$ together with relation (R2a) in all dimensions is needed in order
to define comultiplication $\Delta$.

Recall that there is $n$-framed   Hodge-Tate structure  $h[(Q,\alpha);M]$
related to $H^{2n-1}(CP^{2n-1} \backslash Q, M;\Bbb Q)$ and a choice of a
 family of maximal isotropic subspaces  $\alpha$ on $Q$.

\vskip 3mm \noindent
{\bf Proposition 4.3} {\it There is a homomorphism of abelian groups $h :
S_{2n}(\Bbb C) \rightarrow {\cal H}_{n}$  defined on generators as}
$h : [(Q,\alpha);M] \longrightarrow h[(Q,\alpha);M]$

\vskip 3mm \noindent
{\bf Proof}. We have to show that $h$ maps relations to zero. For R1) and
R2) this is clear from the definitions. If $\cap M_{i} \neq \emptyset$
then  $gr^{W}_{0} h[(Q,\alpha);M] = 0$ and so this $n$-framed
Hodge-Tate structure is equivalent to zero.

Let us prove that
\begin{equation} \label{4.2}
\sum_{j = 0}^{2n}(-1)^{j} h[(Q,\alpha),M^{(j)}] = 0
\end{equation}

One has the following map
\begin{equation} \label{4.3}
 H^{2n-1}({\Bbb C}{\Bbb P}^{2n-1} \backslash Q,\cup_{i = 0}^{2n} M_{i};\Bbb Q)
\stackrel {\oplus r_{j}}{\longrightarrow}
\bigoplus H^{2n-1}({\Bbb C}{\Bbb P}^{2n-1} \backslash Q,M^{(j)};\Bbb Q)
\end{equation}

Consider in the left mixed Hodge structure the $n$-framing $(\omega_{Q},0)$
 Then it is easy to see that $\oplus r_{j}$ induces map of this
Hodge-Tate structure just to the left-hand side of (\ref{4.2})
Therefore (\ref{4.3}) holds. Proposition (4.3) is proved.

Set $S_{0}(\Bbb C) = \Bbb Z$ and $S(\Bbb C)_{\bullet}: =
\oplus_{n=0}^{\infty} S_{2n}(\Bbb C)$

Our next goal is to define on $S(\Bbb C)_{\bullet}$ a structure of graded
commutative Hopf algebra.

The product of elements  $[(Q,\alpha);M] \in S_{2n}(\Bbb C)$
and $[(Q',\alpha');M'] \in S_{2n'}(\Bbb C)$ is defined as follows.
Let $V$ be a vector space of dimension $2n$. Choose a quadratic
equation $\tilde Q$ of the quadric $Q$ and a maximal isotropic subspace
$\tilde \alpha$ from the family  $\alpha$ . Denote by
$\tilde M_{i}$ the codimension $1$ subspace projected to $M_{i}$. There is
an analogous data in a vector space $V'$ related to  $[(Q',\alpha');M']$
. Then in $V \oplus V'$ there are the quadratic form $\tilde
Q \oplus \tilde Q'$, the coordinate simplex consisting of subspaces $\tilde
M_{i} \oplus V'$ and $V \oplus \tilde M_{\bar i}'$  and the maximal
isotropic subspace $\tilde\alpha \oplus \tilde \alpha'$.
 The projectivisation of this data defines a pair  $[(Q
\ast Q';\alpha \ast \alpha'), M \ast M']$   .Its projective equivalence class
does not depend on the choice of $\tilde Q, \tilde Q'$  and $ \tilde \alpha,
\tilde \alpha'$. By definition the product
$$
 [(Q,\alpha);M] \ast  [(Q',\alpha');M']\in S_{2(n+n')}(\Bbb C)
$$
is the class of this pair.

\vskip 3mm \noindent
{\bf Proposition 4.4}. {\it The product is well defined}.

\vskip 3mm \noindent
{\bf Proof}. Clear from the definitions. For example to check
$$
( \sum _{j=0}^{2n} (-1)^{j} [(Q,\alpha);M^{j}])\ast  [(Q',\alpha');M'] = 0
$$
consider $2(n+n')+1$ hyperplanes in $P(V')$   corresponding to $\{ M_{i} \oplus
V'\} \bigcup \{V \oplus  M_{\bar i}'\}$
where $0 \leq i \leq 2n;$ and $ 0 \leq i' \leq {2n-1}$.
Then on the level of generators the left side of (4.6) is just the additivity
relation written for $(Q
\ast Q';\alpha \ast \alpha')$ and these hyperplanes.

Now let us define the comultiplication
$$
\Delta : S_{2n}(\Bbb C) \longrightarrow \bigoplus _{ k = 0}^{n} S_{2k}(\Bbb C)
 \otimes S_{2(n-k)}(\Bbb C)
$$

Let $[(Q,\alpha);M]$ be a generator of $S_{2n}(\Bbb C)$. If $\cap M_{i} \neq
\emptyset$ then $\Delta  [(Q,\alpha);M] = 0$ by definition. Now let us suppose
that $\cap M_{i} \ = \emptyset$. Choose a data $(\tilde Q,\tilde \alpha);
\{\tilde M_{i}\}$   in a vector space $V$ corresponding to
$[(Q,\alpha);M]$ . Let $\tilde Q_{I}$ be the restriction of the quadratic form
$\tilde Q$ to the subspace $\tilde M_{I}$ . Set $\tilde \alpha_{I} : = \tilde
\alpha \cap M_{I}$

Let $\bar I := \{\bar i_{1} < \ldots < \bar i_{2(n-k)}\}$    be the complement
to $I$ in $\{0,...,2n-1\}$. Then there is the
family of hyperplanes $M_{\bar i_{1}} \cap M_{I}, \ldots ,M_{\bar i_{2(n-k)}}
\cap M_{I}$ in $M_{I}$  that will be reffered to as $M_{\{\bar I\}}$  .

Further, let $\tilde Q^{I}$ (resp. $\tilde \alpha^{I}$)   be the quadratic
form ( resp. maximal isotropic subspace) induced in $V/M_{I}$  by
$\tilde Q$ ( resp. $\tilde \alpha$ ). The hyperplanes $M_{i}$ where
$i \in I$ induce the family $M^{I}$
 of hyperplanes in  $P(V/\tilde M_{I}$ . Set $\Delta = \oplus_{k = 0}^{n}
\Delta_{2(n-k),2k}$ and
$$
 \Delta_{2(n-k),2k} ([Q,\alpha];M) : = \sum_{|I| = 2k}([(\tilde Q^{I},\tilde
\alpha^{I});\tilde M^{I}]) \otimes [(\tilde Q_{I},\tilde \alpha_{I});\tilde
M_{\bar I}) \in S_{2k}(\Bbb C) \otimes S_{2(n-k)}(\Bbb C)
$$

There is also a duality $i : S_{2n}(\Bbb C) \rightarrow S_{2n}(\Bbb C)$ defined
as
follows. A nondegenerate quadric $Q \in {\Bbb C}{\Bbb P}^{m}$ provides an
isomorphism
$i_{Q} : {\Bbb C}{\Bbb P}^{m} \rightarrow \hat {\Bbb C}{\Bbb P^{m}}$   .
Namely, let $ x \in {\Bbb C}{\Bbb P}^{m}$. Consider the set of lines through
$x$ tangent to $Q$. The tangency locus is a section of $Q$ by a hyperplane
$i_{Q}(x)$. Let
$v_{j} = \cap_{i \neq j} M_{i}$  be the vertices of simplex $M$. Set
$\hat M : = (i_{Q}(v_{0}), ... , i_{q}(v_{m}))$. Then
$$
i : [(Q,\alpha);M] \longrightarrow [(Q,\beta);\hat M ]
$$

\vskip 3mm \noindent
{\bf Theorem 4.5} a){ \it The comultiplication $\Delta$ is well defined.

b) The multiplication $\ast$ and the comultiplication  $\Delta$ provide the
structure of graded commutative Hopf algebra on $S(\Bbb C)$.

c)The duality $i$ is an antipode, i.e. $i(a \ast b) = i(a) \ast i(b);
\Delta (i(a)) =
t(i \otimes i) \Delta(a)$ ,
where $t$ interchange factors in the tensor product.}

\vskip 3mm \noindent
 {\bf Proof}. a) We have to check that $ \Delta_{2k},_{2(n-k)}$ takes relations
R0)-R3) to zero. This is obvious for R1)-R2) and true by definition for R0). It
remains to check R3).

It is convenient to extend the definition of $\Delta$ to generators
$[(Q,\alpha);M]$ with $\bigcap M_{i} \not= \emptyset$ using (4.8)
where the summation is over all subsets $I$ such that $|I| = 2k$,
codim$M_{I} = 2k$ and $M_{I^{\prime}}\not= M_{I}$ for any $I^{\prime}\supset I$
with $I^{\prime}\not=
I$. It is obvious that then $\Delta[(Q,\alpha);M]=0$ because the second factor
in (4.8) will be always $0$ thanks to R0).

Now let $M_{0},...,M_{2n}$ be hyperplanes in $ {\Bbb C}{\Bbb P}^{2n-1}$ such
that
$Q \bigcap M_{J}$ is a nondegenerate quadric for any subset $J$.

We have to prove that
$$
\Delta_{2k,2(n-k)}(\sum_{j=0}^{2n}(-1)^{j}[(Q,\alpha);M^{(j)}]) = 0
$$
Compute the left side of (4.9) using definition (4.8).The formula we get
is a sum of contributions corresponding to codimension $2k$ edges of the
configuration $(M_{0},...M_{2n})$. Consider one of the contributions related to
a certain
edge $E$.
Let $I$ be the set of all indices $i$ such that $E \subset M_{i}$. If $|I| \geq
2k+2$ then all simplices $M^{(j)}$ are degenerate and so R3) is provided by
R0).
Suppose that $I = 2k+1$. Then in the part of the formula
under consideration the second factor in all summands is the same, and the sum
of the first factors is 0 thanks
to relation R3). Similary if $I =2k$ the first factor is the same, and
the sum of the second ones is $0$ according to R3)

b)We have to prove that $\Delta(a \ast b) = \Delta(a) \ast \Delta(b)$ This is
an immediate consequence of definitions and the following

\vskip 3mm \noindent
{\bf Lemma 4.6}.{\it Let $V = V_{1} \oplus V_{2}, \tilde Q = \tilde
Q_{1} \oplus \tilde Q_{2}, \tilde M =
\tilde M_{1} \oplus \tilde M_{2}$. Suppose that $dimV_{i}$ is odd for $i =
1,2$. Then
the element of $S(C)_{\bullet}$ corresponding to $([Q,\alpha];M)$ is equal to
zero modulo 2-torsion.}

\vskip 3mm \noindent
{\bf Proof}. The transformation $g = (-id_{V_{1}}, id_{V_{2}})$ has determinant
-1
and so interchange connected components of the manifold of maximal isotropic
subspaces of the quadric $Q$ and the simplex $M$. Therefore
$$
([Q,\alpha];M) = ([gQ,g \alpha]; gM) = ([Q, \beta];M) = -([Q, \alpha];M)
$$
Lemma 4.6 is proved.

The coassociativity of $\Delta$ also follows immediately from the definitions
and lemma 4.4.

c) Clear from the definitions.

\vskip 3mm \noindent
{\bf Remark 4.7} (compare with remark 4.2). If we omit from the definition of
group $S_{2n}(\Bbb C)$ orientation of the quadric $Q$ and hence relation (R2a)
we would missed lemma 4.5 and therefore would get $\Delta(a \ast b) \neq
\Delta(a) \ast \Delta(b)$

The comultiplication $\Delta$ can be defined in terms of projective geometry.
Namely $Q_{I}$ and $M_{\{\bar I\}}$ are the quadric and the simplex in $M_{I}$.
To define the second
factor in the wright-hand side of (4.8) consider the projective space formed by
hyperplanes containing $M_{I}$. The tangency condition with $Q$ defines a
quadric in this projective space. Hyperplanes $M_{i}, i\in I$ can be considered
as vertices of a simplex in it. Finally an orientation of the quadric $Q$
provides orientations of both quadrics we get.

Accoding to proposition 4.3 we have the morphism of abelian groups
$$
 h: S(C)_{\bullet} \rightarrow \cal H_{\bullet}
$$

\vskip 3mm \noindent
{\bf Theorem 4.8}. $h$ {\it is a homomorphism of Hopf algebras}.

\vskip 3mm \noindent
{\bf Proof}. Let us prove that $h$ commute with comultiplication. We need the
following elementary facts about the cohomology of quadrics.

a) If $Q \in {\Bbb C}{\Bbb P}^{m}$ is a nonsingular quadric then $H^{i}({\Bbb
C}{\Bbb P}^{m} \backslash Q;\Bbb Q) = 0$
if $m$ is even, $i > 0$ or $m$ is odd, $i \neq 2n-1$ and $H^{2n-1}({\Bbb
C}{\Bbb P}^{2n-1}
\backslash Q; \Bbb Q) = \Bbb Q(-n)$

b) If $Q_{0}$ is a singular quadric in $\Bbb C^{m}$ then $H^{i}(\Bbb
C^{m},Q_{0}; \Bbb Q) = 0$ for all $i$

 Let us compute
\begin{equation} \label{4.4}
H^{\ast}({\Bbb C}{\Bbb P}^{2n-1} \backslash Q, M\backslash (M \cap  Q);\Bbb Q)
\end{equation}

Notice that $(Q,M)$ is a normal crossing divisor. Recall that for a subset $I
\in \{0,...,2n-1\}$ we set $Q_{I} : = M_{I}\cap  Q$. Consider the corresponding
simplicial scheme
$$
{\Bbb C}{\Bbb P}^{2n-1} \backslash Q
\hbox{{$\cdots$}\kern-12.6pt\lower4pt\hbox{$\leftarrow$}
\kern-13.5pt\raise4pt\hbox{$\leftarrow$}} \bigcup_{|I| = 1}  M_{I} \backslash
Q_{I}
\hbox{{$\cdots$}\kern-12.6pt\lower4pt\hbox{$\leftarrow$}
\kern-13.5pt\raise4pt\hbox{$\leftarrow$}}
..
\hbox{{$\cdots$}\kern-12.6pt\lower4pt\hbox{$\leftarrow$}
\kern-13.5pt\raise4pt\hbox{$\leftarrow$}}
\bigcup_{|I| = 2n-1}  M_{I} \backslash Q_{I}
$$

It produces the spectral sequence with
$$
E_{1}^{p,q} = H^{p}(M_{I} \backslash Q_{I}; \Bbb Q)
$$
where $|I| =q$.
This spectral sequence degenerates at $E_{1}$ (because of the weight
considerations). The filtration it induces on (\ref{4.4}) coincides with the
weight filtration. Therefore
$$
gr_{2(n-k)}^{W}H^{2n-1}({\Bbb C}{\Bbb P}^{2n-1} \backslash Q,M;\Bbb Q) =
\bigoplus _{|I| = 2k} H^{2(n-k)-1}(M_{I} \backslash Q_{I};\Bbb Q)
$$

Moreover, for each subset $I$ with $|I| = 2k$ there is an obvious injective
morphism of mixed Hodge structures
$$
i_{I} : h(Q_{I},M_{I}) \rightarrow h(Q,M)
$$

Recall that there are canonical up to a sign vectors
$[\omega_{Q_{I}}]\in gr^{W}_{2(n-k)}h(Q_{I},M_{I})$. If we choose signes,
their images form a basis in $gr^{W}_{2(n-k)}h(Q,M)$. Further, there are
morphisms of mixed Hodge structures
$$
H^{2n-1}({\Bbb C}{\Bbb P}^{2n-1} \backslash Q,M;\Bbb Q) \longrightarrow
H^{2n-1}({\Bbb C}{\Bbb P}^{2n-1} \backslash Q,\cup_{i \in I} M_{i};\Bbb Q)
$$

\vskip 3mm \noindent
{\bf Proposition 4.9}. {\it One has canonical isomorphism of mixed Hodge
structures}

$$
H^{2n-1}({\Bbb C}{\Bbb P}^{2n-1} \backslash Q,\cap_{i \in I} M_{i};\Bbb Q)(n-k)
= h(M^{I},Q^{I})
$$

\vskip 3mm \noindent
{\bf Proof}. By Poincare duality
$$
H^{2n-1}({\bf C}P^{2n-1} \backslash Q,\cup_{i \in I} M_{i};\Bbb Q) =
H^{2n-1}({\Bbb C}{\Bbb P}^{2n-1} \backslash \cup_{i \in I}M_{i}, Q ; \Bbb
Q)^{\ast} (-(2n-1))
$$

Notice that
$$
H^{2n-1}({\Bbb C}{\Bbb P}^{2n-1} \backslash Q,\cup_{i \in I} M_{i} ; \Bbb Q) =
R^{2n-1}f_{\ast}(j_{M \ast}j_{Q!} \Bbb Q)
$$

where
$$
{\Bbb C}{\Bbb P}^{2n-1} \backslash (\cup_{i \in I} M_{i} \cup Q) \stackrel
{j_{Q}}
{\hookrightarrow}    {\Bbb C}{\Bbb P}^{2n-1} \backslash (\cap_{i \in I} M_{i})
\stackrel {j_{M}}
{\hookrightarrow}
 {\Bbb C}{\Bbb P}^{2n-1} \backslash M_{I}
$$

and $f$ is a composition
\begin{equation} \label{4.5}
 {\Bbb C}{\Bbb P}^{2n-1} \backslash M_{I} \stackrel {f_{1}}
{\longrightarrow} P(V/\tilde M_{I}) \backslash M^{I} \stackrel {f_{2}}
{\longrightarrow}
\ast
\end{equation}

Let us compute $R^{2n-1}f_{1\ast}(j_{M \ast}j_{Q!} \Bbb Q)$.
A point $x$ of the plane $P(V/\tilde M_{I})$ corresponds to a plane $H_{x}$ of
dimension dim$M_{I}+1$ containing $M_{I}$. Set $h_{x}:=
H_{x}\backslash M_{I}$. The fiber
$$
i_{x}^{\ast} R^{i}f_{1\ast}(j_{M \ast}j_{Q!} \Bbb Q)
$$
is isomorphic to $H^{i}(h_{x},h_{x} \cap Q)$. Recall that $x$ belongs to the
quadric if and only if $h_{x} \cap Q$ is singular. So the fiber (\ref{4.5}) at
such
points is zero according to b). The fiber (\ref{4.5}) at points
$x \in P(V/\tilde M_{I}) \backslash \Bbb Q^{I} \cap M^{\{I\}}$ is $Q(-k)$ for
$i=2k$ and zero for other $i$. Therefore
$$
Rf_{1\ast}(j_{M \ast}j_{Q!} \Bbb Q) = j_{M^{I} \ast}j_{Q^{I}!} \Bbb Q(-k)[-2k]
$$
and applying again Poincare duality we get proposition 4.9.

\vskip 3mm \noindent
{\bf 8. An analog of Schl\"afli formula.}
To write it I need some facts about geometry of even-dimensional smooth quadric
$Q$.
There are exactly 2 different families $\alpha$ and $\beta$ of maximal
isotropic subspaces on $Q$ (they are of middle dimension). For each isotropic
subspace $\gamma$ of dimension one less there are just  2 maximal isotropic
subspaces $\alpha_{(\gamma)}$ and $\beta_{(\gamma)}$ containing it. They belong
to  families $\alpha$ and $\beta$ respectively. Let $H$ be a generic
codimension 2 hyperplane. Choose an isotropic subspace $\gamma \subset H \cap
Q$. Then there are just two hyperplanes
$H_{\gamma}^{\alpha}$ and $H_{\gamma}^{\beta}$ containing $H$ and tangent to
quadric $Q$ (here $\alpha_{(\gamma}) \subset H_{\gamma}^{\alpha}$).

Let us return to our data $(Q,M)$.Set $M_{ij} = M_i \cap
M_j$. Choose a maximal isotropic subspace $\gamma_{ij} \subset
M_{ij} \cap Q$. Let
$H_{\gamma_{ij}}^{\alpha}$ and $H_{\gamma_{ij}}^{\beta}$ be the
corresponding hyperplanes tangent to $Q$ and containing $M_{ij}$.
The hyperplanes containing $M_{ij}$ form a projective line. Let
$r(M_i, M_j, H_{\gamma_{ij}}^{\alpha}, H_{\gamma_{ij}}^{\beta})$
be the cross -  ratio of 4 points on it corresponding to the hyperplanes

Recall that $v_{\alpha}(Q,M) := \int_{\Delta_M}\omega_Q$ is the
period integral related to
the pair $(Q,M)$ where the sighn
of the form $\omega_Q$ corresponds to the family $\alpha$ of
isotropic planes on $Q$.
Set $Q_{ij} := Q \cap M_{ij}$. The hyperplanes $M_k, k \not =
i,j,$ cut a simplex $M^{ij}$ in $M_{ij}$.

\vskip 3mm \noindent
{\bf Theorem 4.10}
$$
dv_{\alpha}(Q,M) = \sum_{i<j} v_{\gamma_{ij}}(Q_{ij},M^{ij}) \cdot d\log r(M_i,
M_j, H_{\gamma_{ij}}^{\alpha}, H_{\gamma_{ij}}^{\beta})
$$

Notice that if $\gamma'_{ij}$ is a maximal isotropic subspace
from a different family then both the cross-ratio and period $v$
change sign.

Using this formula one can easyly prove that the integrals
$v_{\alpha}(Q,M)$ can
be expressed by the hyperlogarithms of order $n$. Using the fact
that hyperlogarithms of order 3 can be expressed by the classical
trilogarithm (see [G4]) we come to the expression of integrals
$v_{\alpha}(Q,M)$ by means of classical trilogarithms and
products of dilogarithm and logarithms. See [K1] for another
method of proving this result.


\vskip 3mm \noindent
{\bf 9.The category of mixed Tate motives over {\Bbb C}.} According to [B1],
[BD] it
should be a tensor ${\Bbb Q}$-category with a fixed invertible
object $\Bbb Q(1)_{\cal M}$ such that any simple object is isomorphic to
$\Bbb Q(n)_{\cal M}: = \Bbb Q(1)_{\cal M}^{\otimes m} \qquad m\in \Bbb Z$
and these objects are mutually nonisomorphic.

\begin{equation}
Ext^{1}_{\cal M}(\Bbb Q(0)_{\cal M},\Bbb Q(m)_{\cal M}) = 0 \quad \mbox{for}
\quad m \leq 0.
\end{equation}

\begin{equation}
Ext^{i}_{\cal M}(\Bbb Q(0)_{\cal M},\Bbb Q(n)_{\cal M}) =
gr^{\gamma}_{n}K_{2n-i}(\Bbb C) \otimes \Bbb Q
\end{equation}

Furhter, any object in this category carries a canonical increasing finite
filtration (weight filtration), whose graded quotient of degree 2n  is
isomorphic to a direct sum of $\Bbb Q(-n)_{\cal M}$'s.
Finally,
\begin{equation} \label{chert}
Ext^{i}_{\cal M}(\Bbb Q(0)_{\cal M},\Bbb Q(n)_{\cal M}) =
gr^{\gamma}_{n}K_{2n-i}(\Bbb C) \otimes \Bbb Q
\end{equation}

In complete analogy with the case of Hodge-Tate structures one can define
$n$-framed mixed Tate motives. Their equivalence classes form an abelian group
denoted
${\cal A}_{n}$ and ${\cal A}_{\bullet} : = \oplus_{n \geq 0} {\cal A}_{n}$ is
equipped with a structure of a commutative graded Hopf algebra. The category of
mixed Tate motives is canonically equivalent to the category of
finite-dimensional graded
${\cal A}_{\bullet}$-modules. This equivalence sends $\Bbb Q(n)_{\cal M}$
 to the trivial one-dimensional ${\cal A}_{\bullet}$-comodule siting in degree
n and one has
\begin{equation}
Ext^{i}_{\cal M}(\Bbb Q(0)_{\cal M},\Bbb Q(n)_{\cal M}) = H^{i}_{(n)}({\cal
A}_{\bullet},\Bbb Q)
\end{equation}
Here the group in the right side is the degree $n$ part of the
graded vector space $H^{i}({\cal A}_{\bullet},\Bbb Q)$. Therefore
according to axiom (\ref{chert}  ) one should have
\begin{equation}
H^{i}_{(n)}({\cal A}_{\bullet},\Bbb Q) =
gr^{\gamma}_{n}K_{2n-i}(\Bbb C) \otimes \Bbb Q
\end{equation}

Each generator of $S_{2n}(\Bbb C)$ defines an $n$-framed Hodge-Tate structure
realised in cohomology of a (simplicial) algebraic variety and therefore an
$n$-framed mixed Tate motive.
Similar arguments provides us with canonical homomorphism of
graded Hopf algebras
\begin{equation}  \label{4.6}
S(\Bbb C)_{\bullet} \longrightarrow {\cal A}_{\bullet}
\end{equation}

Therefore should be canonical homomorphisms
\begin{equation} \label{4.7}
H^{i}_{(n)}(S(\Bbb C)_{\bullet},\Bbb Q) \longrightarrow
gr^{\gamma}_{n}K_{2n-i}(\Bbb C) \otimes \Bbb Q
\end{equation}

One can compute the left side as the cohomology of degree $n$ part of the cobar
complex $S^{\bullet}(n)$ for the graded Hopf algebra $S(\Bbb C)_{\bullet}$:
\begin{equation}
S_{n}(\Bbb C) \longrightarrow \oplus_{i_{1} + i_{2}= n} S_{i_{1}}(\Bbb C)
\otimes S_{i_{2}}(\Bbb C) \longrightarrow \oplus_{i_{1} + i_{2} + i_{3}= n}
S_{i_{1}}(\Bbb C)_{i_{1}}
\otimes S_{i_{2}}(\Bbb C) \otimes S_{i_{3}}(\Bbb C) \longrightarrow ...
\end{equation}

Here the left group is placed at degree 1 and the coboudary is of degree +1.
Therefore we arrive to the following

\vskip 3mm \noindent
{\bf Conjecture 4.11 }. a){\it There is canonical homomorphism}
\begin{equation}  \label{ss}
H^{i}(S^{\bullet }(n)) \otimes \Bbb Q \longrightarrow
gr^{\gamma}_{n}K_{2n-i}(\Bbb C)
\otimes \Bbb Q
\end{equation}
b) {\it It is injective for $i= 1$}

\vskip 3mm \noindent
{\bf Problem 4.12 } {\it Wherether it is true that homomorphism (\ref{ss}) is
an isomorphism?}

Notice that (\ref{4.7}) is an isomorphism if and only if (\ref{4.6})
is an isomorphism. I believe that (\ref{4.6}) is at least injective; this
provides the part b) of conjecture 4.11.

\vskip 3mm \noindent
{\bf Theorem 4.13 } {\it Conjecture 4.11 valid for $n$ = 2}

Let me deduce conjectures 2.4 and 2.7 from conjecture 4.11. Let us represent
${\Bbb C}{\Bbb P}^{n}$ as a complexification of ${\bf RP}^{n}$. Then the
complex conjugation acts on ${\Bbb C}{\Bbb P}^{n}$. Set


\vskip 3mm \noindent
{\bf Theorem 4.14}. a){\it There are canonical homomorphisms
\begin{equation}  \label{p1}
{\cal P}({\cal H}^{2n-1}) \stackrel {\psi_{{\cal H}}(n)}{\longrightarrow}
(S_{2n}(\Bbb C) \otimes \Bbb Z(n))^{-}
\end{equation}
\begin{equation}  \label{p2}
{\cal P}(S^{2n-1})  \stackrel {\psi_{S}(n)}{\longrightarrow} (S_{2n}(\Bbb C)
\otimes \Bbb Z(n))^{+}
\end{equation}

b) They transform Dehn invariant just to the comultiplication, providing the
following  commutative diagram
in hyperbolic case

\vskip 3mm
\begin{center}
\begin{picture}(100,100)(25,-10)
\put(0,75){$
{\cal P}({\cal H}^{2n-1})
$}
\put(50,80){\vector(1,0){90}}
\put(0,0){$
S_{2n}(\Bbb C)
$}
\put(150,75){$
\oplus {\cal P}({\cal H}^{2k-1}) \otimes {\cal P}(S^{2(n-k)-1})
$}
\put(50,0){\vector(1,0){90}}
\put(160,0){$
\oplus S_{2k}(\Bbb C) \otimes S_{2(n-k)}(\Bbb C)
$}
\put(25,65){\vector(0,-1){45}}
\put(180,65){\vector(0,-1){45}}
\put(-10,40){$
 {\psi}_{{\cal H}}(n)
$}
\put(185,40){$ {\psi}_{{\cal H}}(k) \otimes {\psi}_{S}(n-k)$}
\put(90,5){$\Delta$}
\put(87,90){$D_{2n-1}^{H}$}
\end{picture}
\end{center}
and a similar diagram in the spherical case.}

c){\it The homomorphisms (\ref{p1}) , (\ref{p2}) provide the
  following homomorphisms from the hyperbolical  and spherical
  Dehn complexes}  (see (\ref{dehn1})
\begin{equation}  \label{d1}
\psi_{{\cal H}}^{\bullet}(n): {\cal P}_{{\cal H}} ^{\bullet}(n)
\longrightarrow (S^{\bullet}(n) \otimes \Bbb Z(n))^{-}
\end{equation}
\begin{equation}
\psi_{S}^{\bullet}(n): {\cal P}_{S}^{\bullet}(n) \longrightarrow
(S^{\bullet}(n) \otimes \Bbb Z(n))^{+}
\end{equation}

\vskip 3mm \noindent
{\bf Proof}. Follows from the definition of $\Delta$ and theorem 4.5a) The
action of the complex conjugation on the image of scissor congruence groups
computed easyly looking on the frames.
Notice that homomorphisms $\psi_{{\cal H}}(n)$ and $\psi_{S}(n)$  are injective
but certainly not isomorphisms unless for $\psi_{{\cal H}}(n)$

Combining this homomorphisms with conjectures 4.11a) and 4.11b)
we get conjecture 2.7 and conjecture 2.4 respectively.

\vskip 3mm \noindent
{\bf 10. The structure of the complexes $S^{\bullet}(n)$ for $n
  \leq 3$}.

{\bf Conjecture 4.15}. {\it There exists canonical homomorphisms of complexes
$$
\begin{array}{ccc}
S_2(\Bbb C)& \longrightarrow & \Bbb C^* \otimes \Bbb C^* \\
&&\\
\downarrow s_2&& \downarrow \\
&&\\
B_2(\Bbb C)& \longrightarrow & \Lambda^2 \Bbb C^*
\end{array}
$$
and
$$
\begin{array}{ccccc}
S_3(\Bbb C)& \longrightarrow & \Bbb C^*\otimes  S_2(\Bbb C)
\oplus  S_2(\Bbb C) \otimes \Bbb C^* & \longrightarrow &
\Bbb C^*\otimes \Bbb C^*\otimes \Bbb C^*\\
&&&&\\
\downarrow s_3&& \downarrow id \otimes s_2 + s_2 \otimes id && \downarrow\\
&&&&\\
B_3(\Bbb C)& \longrightarrow & B_2(\Bbb C) \otimes \Bbb C^* &
\longrightarrow&\Lambda^3 \Bbb C^*
\end{array}
$$
which induce isomorphisms on homology modulo torsion}.

Here the group $B_2(\Bbb C)$ was introduced in 1.4 and the group
$B_3(\Bbb C)$ is an
explicit version of the group ${\cal B}_3(\Bbb C)$ from
s.1.5, see [G2]. The differentials in the $B$-complexes are defined by the
formulas $\{x\}_2 \longmapsto (1-x) \wedge x$ and $\{x\}_3
\longmapsto \{x\}_2 \otimes x$, $\{x\}_2\otimes  y \longmapsto
(1-x) \wedge x \wedge y$.
A similar results, of cource, should be valid for any algebraicly
closed field.

The most interesting problem here is a construction of an
{\it explicit} homomorphisms  $S_3(\Bbb C) \longrightarrow B_3(\Bbb C)$
and $S_2(\Bbb C) \longrightarrow B_2(\Bbb C)$ which commutes with
the differentials. I believe that these homomorphisms should have
a very beautyfull geometrical description. For example,
the homomorphism $S_2(\Bbb C) \longrightarrow B_2(\Bbb C)$ is
easy to define for simplices whose vertices are on the quadric
$Q$. Namely, $Q =\Bbb CP^1 \times \Bbb CP^1$, so such a simplex
is defined by a 4-tuple
$[(x_1,y_1),(x_2,y_2),(x_3,y_3),(x_4,y_4)]$ of points of $P^1
\times P^1$. Its
image in $B_2(\Bbb C)$ should be
$\{r(x_1,x_2,x_3,x_4)\}_2 - \{r(x_1,x_2,x_3,x_4)\}_2$.

The homomorphisms $s_n$ ($n=2,3$) will lead to explicit
formulas for volumes of noneuclidean simplices via classical di
and trilogarithms (compare with the work of R.Kellerhals
[K1]). More precisely,
the composition
$$
{\cal P}({\cal H}^{2n-1}) \longrightarrow S_n(\Bbb C)
\stackrel{s_n}{\longrightarrow} B_n(\Bbb C)
\stackrel{{\cal L}_n}{\longrightarrow} \Bbb R
$$
 should coinside with the volume homomorphism.
The kernel of the maps $s_2, s_3$ should be equal to $S_1(\Bbb C)
 \ast S_1(\Bbb C)$ and $S_1(\Bbb C)
 \ast S_2(\Bbb C)$.

\vskip 3mm \noindent
{\bf 11. A construction of an element in $Ext^{1}_{\cal M}(\Bbb Q(0)_{\cal
M},\Bbb Q(n)_{\cal M})$   corresponding to a hyperbolic $(2n-1)$- manifold
.}
Recall that a hyperbolic manifold $M^{2n-1}$ produces an element
$s(M^{2n-1}) \in {\cal P}({\cal H}^{2n-1})$  that is
obtained by a decomposition of the manifold on geodesic simplices (see s.2.1).

Proposition 2.3 provides $s(M^{2n-1}) \in Ker D^{H}_{2n-1}$ and so
thanks to theorem 4.14b) its
 image $\psi_{n}( s(M^{2n-1}))$ under homomorphism (\ref{d1}) belongs to
Ker$\Delta$.

Therefore according to theorem 4.8 the comultiplication $\nu$ of the
corresponding n-framed Hodge-Tate structure
$h \circ \psi_{{\cal H}}(n) ( s(M^{2n-1}))$ is also $0$.

By the definition of the graded Hopf algebra $\cal H_{\bullet}$ one has
$h \circ \psi_{{\cal H}}(n)( s(M^{2n-1})) \in {\cal H}_n$. Look at the cobar
complex for the Hopf algebra $\cal H_{\bullet}$:
\begin{equation}
\cal H_{\bullet} \stackrel {\nu}{\longrightarrow} {\cal H}_{\bullet}
\otimes\cal H_{\bullet} \longrightarrow ...
\end{equation}

The fact that the comultiplication $\nu$ of $h \circ \psi_{{\cal H}}(n) (
s(M^{2n-1}))$ is zero  just means that
$h \circ \psi_{{\cal H}}(n) ( s(M^{2n-1}))$ is a 1-cocycle in this complex.

Therefore according to theorem 4.1  the cohomology class of this 1-cocycle
 provides us an element of $Ext^{1}_{{\cal H}_{T}}(\Bbb Q(0)_{{\cal H}_{T}},
\Bbb Q(n)_{{\cal H}_{T}})$ that
is by construction of algebraic-geometrical origin. So we have
 constructed a  motivic extension promised in s. 1.5.

\vskip 3mm \noindent
{\bf  12. Problem 2.7 has positive answer for $n=2$.}

\vskip 3mm \noindent
{\bf Theorem 4.16}([D],[DS1],[DPS]) {\it The sequences}
\begin{equation} \label{h}
0 \longrightarrow H_{3}(SL_{2}(\Bbb C))^{-} \longrightarrow
{\cal P}({\cal H}^{3}) \stackrel {D^{H}_{3}}{\longrightarrow }
\Bbb R \otimes S^{1}
\longrightarrow H_{2}(SL_{2}(\Bbb C))^{-} \longrightarrow 0
\end{equation}
\begin{equation}
0 \longrightarrow H_{3}(SU(2))\longrightarrow
{\cal P}(S^{3}) \stackrel {D^{S}_{3}}{\longrightarrow }
S^{1} \otimes S^{1}
\longrightarrow H_{2}(SU(2)) \longrightarrow 0
\end{equation}
{\it are exact}.

According to [Su 2] and [Sah 3] $H_{3}SL_{2}(\Bbb C) =
  K_{3}^{ind}(\Bbb C)$ modulo torsion. Further, $H_{2}SL_{2}(\Bbb C)
= K_{2}(\Bbb C), \qquad gr^{\gamma}_{2}K_{2}(\Bbb C) = K_{2}(\Bbb C)$ and
$gr^{\gamma}_{2}K_{3}(\Bbb C) = K_{3}^{ind}(\Bbb C)$.
Therefore the complex
\begin{equation} \label{dul}
{\cal P}({\cal H}^{3}) \stackrel {D^{H}_{3}}{\longrightarrow }
\Bbb R \otimes S^{1}
\end{equation}
computes $K_{3}^{ind}(\Bbb C)^{-}$ modulo torsion and $K_{2}(\Bbb C)^{-}$.

Moreover, the complex (\ref{dul}) is just the ``-'' part of the following Bloch
complex [D]:
\begin{equation} \label{bl}
0 \longrightarrow B_{2}(\Bbb C)\stackrel {\delta}{\longrightarrow }
\wedge^{2} \Bbb C^{*}\longrightarrow 0   \qquad \delta \{x\}_2
:= (1-x) \wedge x
\end{equation}

(By Matsumoto theorem Coker$\delta = K_{2}(\Bbb C)  $ and according to [Su 2]
and [Sah 3] Ker$\delta =  K_{3}^{ind}(\Bbb C)$ modulo torsion.)

The construction of the map of complexes is provided by the
isomorphisms ${\cal P}({\cal H}^{3})
= {\cal P}(\partial {\cal H}^{3})$ and $\partial {\cal H}^{3} = {\Bbb C}{\Bbb
P}^{1}$;
 the coincidence of differentials is not obvious and follows from
calculations.

Finally, in [DPS] was given a rather involved construction (using  the Hopf map
$S^{3} \longrightarrow S^{2}$) of a homomorphism of the spherical complex
\begin{equation}
{\cal P}(S^{3}) \stackrel {D^{S}_{3}}{\longrightarrow }
S^{1} \otimes S^{1}
\end{equation}
to the ``+'' part of
(\ref{bl}). It is not an isomorphism but induces isomorphisms on cohomologies
([WS],
[DPS]).

\vskip 3mm \noindent
{\bf 13. A symplectic approach to the Hopf algebra $S(\Bbb C)_{\bullet}$}. Let
$W_{2n}
$ be a symplectic vector space decomposed to a direct sum of Lagrangian
subspaces $W_{2n} = E \oplus F$. Suppose also that coordinate hyperplanes $E
_{1}, ... , E_{n}$ in $E$ are given. The symplectic structure provides an
isomorphism $F = E^{\ast}$ and therefore there are dual hyperplanes $F_{1}, ...
,F_{n}$ in $F$.

The subgroup of $Sp(W_{2n})$ preserving this data is an $n$-dimensional torus
$T_{n}$ ( a maximal Cartan subgroup).

Let $L_{n}^{0}$ be the manifold of all Lagrangian subspaces
$H \in W_{2n}$ in generic position with respect to coordinate hyperplanes
$E_{i} \oplus F, E \oplus F_{j}$. Each of them can be considered as a graph
of a map $h: E \rightarrow F$ and hence defines a bilinear form
$h \in E^{\ast} \otimes F^{ \ast}$. The condition that $H$ is isotropic
just means that $h$ is symmetric bilinear form.

{\bf Lemma 4.17}. {\it Points of $L_{n}^{0}/T_{n}$ are in 1-1 correspondence
with configurations (i.e. projective equivalence classes) of pairs $(Q;M)$
where $Q$ is a nondegenerate quadric and $M$ is a simplex in generic position
with respect to $Q$ in $P(E)$}.

{\bf Proof}. The coordinate hyperplanes $E_{1}, ... ,E_{n}$ define a simplex in
$P(E)$. The torus $T_{n}$ is the subgroup of all transformations in $GL(E)$
preserving this simplex. Lemma 4.17 follows from these remarks.

\section{ Proof of the theorem 2.5} 

{\bf 1. Some results on the  $t$-structure on a triangulated
  category}.
 Recall that a $t$-structure on a triangulated  category  $\cal D$
is a pair of subcategories ${\cal D}^{\leq 0}$,${\cal D}^{\geq
  1}$ satisfying the following
conditions:

1)$Hom_{\cal D}(X,Y) =0$ for all objects $X \in Ob
{\cal D}^{\leq 0}$ and $Y \in {\cal D}^{\geq 1}$.

2)For any $X \in Ob {\cal D}$ there exists an exact triangle
$$
X_{\leq 0} \longrightarrow X \longrightarrow X_{\geq 1}
\longrightarrow X_{\leq 0}[1]
$$
with $X_{\leq 0} \in {\cal D}^{\leq 0}$ and $X_{\geq 1} \in {\cal
  D}^{\geq 1}$.

3)${\cal D}^{\leq 0} \subset {\cal D}^{\leq 1}$ and
${\cal D}^{\geq 1} \subset {\cal D}^{\leq 0}$

Here ${\cal D}^{\leq a}:= {\cal D}^{\leq 0}[-a]$, ${\cal D}^{\geq
  a}:= {\cal D}^{\geq 0}[-a]$, ${\cal D}^{[a,b]}:= {\cal D}^{\leq
  b}\cap {\cal D}^{\geq a}$ and ${\cal D}^{a}:= {\cal
  D}^{[a,a]}$.

The exact triangle in (2) is defined uniquely up
to an isomorphism and depends functorially on $X$.
The subcategory ${\cal D}^{0}$ is called the heart of a
$t$-structure. It is an abelian category ([BBD]).

Let ${\cal A}$ and ${\cal B}$ be two sets of isomorphism classes
of objects in ${\cal D}$. Denote by
${\cal A} \ast {\cal B}$ the set of all objects $X$ in $Ob {\cal
  D}$ which can be included into an exact triangle $A
\longrightarrow X \longrightarrow  B \longrightarrow
A[1]$ with $A \in {\cal A}, B \in {\cal B}$.

\begin{lemma} \label{asso}
$({\cal A} \ast {\cal B}) \ast {\cal C} = {\cal A} \ast ({\cal B}
\ast {\cal C})$
\end{lemma}

{\bf Proof}. Follows from the octahedron lemma, see  [BBD].

\begin{theorem} \label{trc}
 a)Let $\cal D$ be a triangulated  $\Bbb Q$-category
and ${\cal Q}$ be a  full semisimple subcategory generating ${\cal D}$ as
a triangulated  category. Suppose that for any two objects
$Q_1,Q_2 \in Ob {\cal Q}$ one has
\begin{equation} \label{negext}
Hom_{\cal D}^{-i}(Q_1,Q_2) = 0 \qquad i> 0
\end{equation}
Then there is  canonical t-structure on ${\cal D}$   with the
abelian heart ${\cal M} = \cup{\cal Q} \ast {\cal Q} \ast ... \ast {\cal
  Q}$.

b) If in addition $Hom_{\cal D}^{i}(Q_1,Q_2) = 0$ for
$i \geq 2$, then the
 tensor category ${\cal D}$ is
equivalent to the derived
category of ${\cal M}$.
\end{theorem}

\begin{theorem} \label{trc1}
Let $\cal D$ be a triangulated  tensor $\Bbb Q$-category
and ${\cal Q}$ be a  full semisimple subcategory generated
by non isomorphic
objects $\Bbb Q(m)$, $m \in \Bbb Z$ such that $\Bbb Q(1)$
is invertible,
$\Bbb Q(m) = \Bbb Q(1)^{\otimes m}$, and
$$
Hom_{\cal D}^{i}( \Bbb Q(m),\Bbb Q(n)) = 0 \quad \mbox{if} \quad  m > n, \qquad
Hom_{\cal D}( \Bbb Q(0),\Bbb Q(0)) =  \Bbb Q
$$
Then the  abelian heart ${\cal M}$
from the  theorem (\ref{trc}) is a tensor category. It is  equivalent to
the tensor category of finite dimensional
representations of a certain free negatively graded  (pro)-Lie algebra.
\end{theorem}

Let ${\cal D}$ be a triangulated category and ${\cal M} \subset {\cal D}$
be a full subcategory.

\begin{theorem} \label{trc3}
${\cal M}$ is a heart for the unique bounded t-structure on
${\cal D}$ if and only if

i) ${\cal M}$ generates ${\cal D}$ as a triangulated category

ii) ${\cal M}$ is closed with respect to extensions

iii) $Hom_{\cal D}(X,Y[i]) = 0 \quad \mbox{for any} \quad X,Y \in
{\cal M}, \quad i<0 $

iv) ${\cal M} \ast {\cal M}[1] \subset {\cal M}[1] \ast {\cal M}$
\end{theorem}

For the proof of this theorem see [BBD] or [P]. Let me scetch
the construction of the $t$-structure on ${\cal D}$.  We will use
the following

\begin{lemma} \label{dirs}
Let
$$
X \longrightarrow Y \longrightarrow Z
\stackrel{f}{\longrightarrow} X[1]
$$
be an exact triangle. Then $f=0$ if and only if $Y =X\oplus Z$
\end{lemma}

{\bf Proof}. Let us show that  $f=0$ implies $Y =X\oplus Z$.
The composition of the identity morphism
$Z \stackrel{id}{\longrightarrow} Z$ with $f$
is zero, so one has  a morphism
$g: Z  \longrightarrow Y$ making the following diagram commutative:
$$
\begin{array}{ccccccc}
 X & \stackrel{a}{\longrightarrow} & Y& \longrightarrow & Z &
\stackrel{f}{\longrightarrow}& X[1]\\
&&&&&&\\
&&&\nwarrow g& \uparrow id& \nearrow =0\\
&&&&&&\\
&&&&Z&&
\end{array}
$$
So there is a morphism $X \oplus Z
\longrightarrow X$. The universality property of this morphism
follows immediately from $f=0$.

Let ${\cal D}^{[a,b]}$ be the minimal full subcategory of ${\cal D}$
containing ${\cal M}[-i]$
for $a \leq i \leq b$ and closed under extensions. The
subcategories ${\cal D}^{\leq 0}$ and ${\cal D}^{\geq 1}$ are
defined similarly and satisfy 1) thanks to iii). Notice that
$Hom_{{\cal D}}({\cal M},{\cal M}[-i]) =0$ implies ${\cal
  M} \ast {\cal M}[i+1] = {\cal M} \oplus {\cal M}[i+1]$ for $
i>0$ by lemma \ref{dirs}. So ${\cal M} \ast {\cal
  M}[n] \subset {\cal M}[n] \ast {\cal M}$ for any $n>0$. Using
ii) and
lemma \ref{asso} we get ${\cal D}^{[a,b]} = {\cal M}[-a]\ast
{\cal M}[-a-1]\ast ... \ast {\cal M}[-b]$. This proves 2). It
remains to show ${\cal D}^{\leq 0} \cap {\cal D}^{\geq 1}
\subset {\cal D}^{[a,b]}$. If $X = Y_c \ast Y_{c+1} \ast ... \ast Y_d$
where $Y_i \in {\cal M}[-i] \in {\cal D}^{\geq a}$ then
$Hom(Y_c,X) =0$ for $c<a$ and so we get from the exact triangle
$$
Y_c \longrightarrow X \longrightarrow  Y_{c+1} \ast ... \ast Y_d
\longrightarrow Y_c[1]
$$
 that $Y_{c+1} \ast ... \ast Y_d = X \oplus Y_c[1]$. Now
 $Hom_{{\cal D}}(Y_c[1],Y_{c+1} \ast ... \ast Y_d) =0$ by iii)
 and so $Y_c
 = 0$.

{\bf Proof of theorem\ref{trc}}. a).
Let us  use  theorem \ref{trc3}.
The properties i)-ii) are obvious. iii) is easy to prove
using (\ref{negext}), so we have to check
only iv). Using the associativity of $\ast$ operation we see that
one has to prove only that  $Q_1 \ast Q_2[1] \subset Q_2[1] \ast Q_1
$ for any two simple objects in ${\cal Q}$. One has
$$
Q_1  \longrightarrow X \longrightarrow Q_2[1]
\stackrel{f}{\longrightarrow} Q_1[1]
$$
There are only two possibilities for
$f$:

1. $f$ is an isomorphism; then $X=0$.

2. $f=0$; then $X = Q_1 \oplus Q_2[1]$ by lemma \ref{dirs}.

 {\bf Part b)}.
\begin{proposition} \label{ex}
a) Let ${\cal D}$ be a triangulated category and ${\cal M}$ be the
heart of a $t$-structure on ${\cal D}$. Suppose that
one has
\begin{equation} \label{ext2}
Hom^i_{{\cal D}}(X,Y)
=0 \quad \mbox{for any} \quad X,Y \in Ob{\cal M}, \quad i >1
\end{equation}
Then the  category
$D^b({\cal M})$ is
 equivalent to the  category ${\cal D}$.

b) If we suppose in addition that the triangulated category ${\cal D}$ is
a subcategory of a derived category $D^b({\cal N})$ for some abelian category
${\cal N}$, then $D^b({\cal M})$ is equivalent to ${\cal D}$ as a
triangulated category.
\end{proposition}

{\bf Proof}. Let $\tilde {\cal D}$ be the full subcategory of
${\cal D}$
whose objects are direct sums $\oplus A_i[-i]$ where $A_i \in
{\cal D}^{-i}$,
 i.e. $H^j_{{\cal M}}(A_i) =0$ for $j \not = -i$. There is a
 canonical functor $i: \tilde {\cal D} \hookrightarrow {\cal D}$.
Then (\ref{ext2}) is a necessary and
sufficient condition for $i$ to be an equivalence of triangulated
categories.
Indeed, let us show that every object $Y$ in ${\cal D}$ is isomorphic
 to $i(X)$ for some $X$. We may suppose $Y \in {\cal D}^{[0,n]}$ and will use
 induction by $n$. Consider the exact triangle
$$
H^0_{{\cal M}}(Y) \stackrel{f}{\longrightarrow} Y \longrightarrow C
\longrightarrow H^0_{{\cal M}}(Y)[1]
$$
provided by the canonical morphism $f: H^0_{{\cal M}}(Y) \longrightarrow Y
$. Then $C \in {\cal D}^{[1,n]}$ and $H^0_{{\cal M}}(Y)[1] \in
{\cal D}^{-1}$,
so it is easy to see that the morphism  $C
\longrightarrow H^0_{{\cal M}}(Y)[1]$ must equal to $0$ because
$Hom^{\geq 2}$ are zero. Therefore $Y = H^0_{{\cal M}}(Y) \oplus C$ by the
lemma
above.

Let us show that $D^b({\cal M})$ has cohomological dimension
one, i.e.
\begin{equation} \label{ext2!}
Hom^i_{D^b({\cal M})}(X,Y)
=0 \quad \mbox{for any} \quad X,Y \in Ob{\cal M}, \quad i >1
\end{equation}
Notice
that for any $X,Y \in {\cal M}$ one has
$$
 Hom_{D^b({\cal M})}^i(X,Y) = Ext_Y^i(X,Y)
$$
where $Ext_Y^i(X,Y)$  is the Yoneda $Ext$-groups in ${\cal
  M}$. Also
$Hom_{{\cal D}}^1(X,Y) = Ext_Y^1(X,Y)$.

One has $Ext^2_Y(X,Y) \subset Ext^2_{{\cal D}}(X,Y)$.
Therefore $Ext^2_{{\cal D}}(X,Y)=0$ implies $Ext^2_Y(X,Y) =0$.
Any
element in $Ext_Y^n$ can be represented as a  product of certain
elements from $Ext_Y^1$ and $ Ext_Y^{n-1}$. So if $Ext^2_Y(X,Y) = 0$ for any 2
objects $X,Y$ in ${\cal
  M}$, then $Ext^i_Y(X,Y) = 0$ for all $i \geq 2$. So we have  an
equivalence of the categories.  Part b) is proved using [BBD].
The proposition is proved.

The proof of the theorem (\ref{trc1}) is rather standard.
One shows that ${\cal M}$ is a mixed Tate category
(see  [BD] or [G1] for the definitions) thanks to the conditions
imposed on $Hom$'s between $\Bbb Q(i)$'s; then the
Tannakian formalism leads to the theorem (\ref{trc1})
(see again [BD] or [G1]).

{\bf 2. An application: the abelian category of mixed Tate
  motives over a number field}.
I will work with the category ${\cal D}{\cal M}_F$  of triangulated
 mixed motives over a field $F$ from [V].  An object in ${\cal D}{\cal M}_F$
is a ``complex'' of regular (but not necessarily projective) varieties $X_1
\longrightarrow X_2
\longrightarrow ... \longrightarrow X_n$ where the morphisms are given by
finite
correspondences and the composition of any two successive morphisms is
zero.

A pair
 $(\Bbb P^n \backslash Q,L)$ where $Q$ is a nondegenerate quadric and $L$
is a simplex provides an object $m(Q,L)$ in
${\cal D}{\cal M}_F$. Namely, for an  subset set $I= \{j_1 < ... <j_{n-i}\}$
of $ \{0,1,...,n \}$ let $L(I):=
L_{j_1} \cap L_{j_2} \cap ... \cap
L_{j_{n-i}}$ be the corresponding $i$-dimensional face of $L$ and
$L^Q(I):=L(I) \backslash (Q \cap L(I))$. Let $L^Q(i):= \cup_{|I| =
n-i}L^Q(I)$. Then $m(Q,L):=  Hom({\tilde m}(Q,L),\Bbb Q)$
 Here $Hom$ is the inner Hom in the category ${\cal
 D}{\cal M}_F$ and
$$
{\tilde m}(Q,L):= L^Q(0)  \longrightarrow L^Q(1) \longrightarrow ...
\longrightarrow L^Q(n-1) \longrightarrow \Bbb P^n \backslash Q
$$
where the first group is sitting in degree $0$ and the
differentials decrease the
degree and given by the usual rules in the simplicial resolution.

There are the objects $\Bbb Q(n) \in {\cal D}{\cal M}_F$ satisfying
almost all
the needed properties including the relation with K-theory:
\begin{equation} \label{vvvvv}
Ext^i_{{\cal D}{\cal M}_F}(\Bbb Q(0), \Bbb Q(n)) =
gr^{\gamma}_nK_{2n-i}(F)\otimes \Bbb Q
\end{equation}
The formula above follows from the results of Bloch, Suslin and
Voevodsky; the key step is the relation of  Higher Chow groups
and algebraic K-theory proved by Bloch (see [Bl3] and  the moving
lemma in [Bl4]). For the relation between the   Higher Chow groups
and motivic cohomology (i.e. the left hand side in (\ref{vvvvv})) see
[V], Proposition 4.2.9 (Higher Chow groups = Borel-Moore motivic homology)
and [V], Proposition 4.3.7 (duality for smooth varieties).

The only serious problem  is
the Beilinson-Soul\'e vanishing conjecture which should guarantee
that the negative $Ext$'s are zero.
However if $F$ is a number field the vanishing conjecture   follows from the
results of Borel and Beilinson ([B2], [Bo2]). So $Hom_{{\cal D}{\cal
    M}}^{-i}(\Bbb Q, \Bbb Q(n))=0$ for all $i,n>0$. The category ${\cal D}{\cal
M}_F$ is a subcategory of the
  derived category of ``sheaves with transfers'', see [V], so we
may apply part b) of the proposition 5.6.
  Therefore we can apply the ``category
machine'' from s. 5.1 and get the abelian category  ${\cal M}_T(F)$
of mixed Tate
  motives over a number field $F$ with {\it all} the needed
  properties. In particulary we have  the Hopf algebra
  ${\cal A}_{\bullet}(F)$ provided by the Tannakian formalism.

For any embedding $\sigma: F \hookrightarrow \Bbb C$
there is the realisation functor  $H_{\sigma}$
from the abelian category of mixed motives over a number field
$F$ to the abelian category of mixed Hodge Tate structures:
$$
H_{\sigma}: {\cal M}_T(F) \longrightarrow {\cal H}_T
$$
It follows from the Borel theorem (injectivity of
the regulator map on $K_{2n-1}(F)\otimes \Bbb Q$ for number fields) that
$\oplus_{\sigma} H_{\sigma}$ induces an injective map on
\begin{equation} \label{BeBo1}
\oplus_{\sigma} H_{\sigma}: Ext_{{\cal M}_T(F)}^1(\Bbb Q(0), \Bbb
Q(n))  \hookrightarrow \oplus_{\sigma}
Ext_{{\cal H}_T(F)}^1(\Bbb Q(0), \Bbb Q(n))
\end{equation}

{\bf 3.  Proof of  the theorem (\ref{Theorem 2.4}): the final step}.
The Hodge realisation provides a  diagram
$$
\begin{array}{ccc} \label{codi1}
S(F)_{\bullet} & \stackrel{\Delta}{\longrightarrow} &   S(F)_{\bullet} \otimes
S(F)_{\bullet}\\
&&\\
\downarrow&&\downarrow\\
&&\\
\oplus_{\sigma} {\cal H}_{\bullet }&
\stackrel{\nu}{\longrightarrow} &
\oplus_{\sigma}{\cal H}_{\bullet} \otimes {\cal H}_{\bullet}
\end{array}
$$
Here the vertical arrows are embeddings because (\ref{BeBo1}) is
injective.  This diagram is commutative thanks to the main
results of the
chapter 4 (especially theorem 4.8).

Further, one has the diagram
$$
\begin{array}{ccc} \label{codi1}
S_{\bullet}(F) & \stackrel{\Delta}{\longrightarrow} &   S_{\bullet}(F) \otimes
S_{\bullet}(F)\\
&&\\
\downarrow&&\downarrow\\
&&\\
{\cal A}_{\bullet}(F)& \stackrel{\Delta}{\longrightarrow} & {\cal
  A}_{\bullet}(F) \otimes {\cal
  A}_{\bullet}(F)\\
&&\\
\downarrow&&\downarrow\\
&&\\
\oplus_{\sigma} {\cal H}_{\bullet }&
\stackrel{\nu}{\longrightarrow} &
\oplus_{\sigma}{\cal H}_{\bullet} \otimes {\cal H}_{\bullet}
\end{array}
$$
where the composition of vertical arrows coincides with the
corresponding
vertical arrow in the previous diagram. (Abusing notations we
denoted   the comultiplication   in the two different Hopf
algebras by the same letter
$\Delta$).
This, together with
injectivity of vertical arrows in the bottom square
of the second diagram implies the
commutativity of the upper square of that diagram. Theorem
(\ref{Theorem 2.4}) follows from the commutativity of the upper square
in the second diagram. Indeed,
the kernel of the middle horisontal arrow coincides with
$K_{2n-1}(F)\otimes \Bbb Q$ and the Beilinson regulator comes
from the bottom square of that diagram.

{\bf Acknoledgements.} This work was essentially done during my stay in
MPI(Bonn)
MSRI(Berkeley) in 1992 and  supported by NSF Grant DMS-9022140.
I am
grateful to
these institutions for hospitality and support. I was also
supported by the NSF Grant  DMS-9500010 in 1995.

I am indebted to M.
Kontsevich, and D.B.Fuchs
for useful discussions and to  L.Positselsky and V.Voevodsky for
helpful suggestions regarding
the last section. I am very gratefull to
A.Borel who read a preliminary version of this paper and pointed out many
misprints and some errors.

\vskip 3mm \noindent
{\bf REFERENCES}
\begin{itemize}
\item[{[Ao]}] Aomoto K.: {\it Analitic structure of Schlafli
    function}, Nagoya Math. J. 68 (1977),
1-16
\item[{[B1]}] Beilinson A.A.: {\it Height pairings between
algebraic cycles}, Lecture Notes in Math. N. 1289, (1987), p.
1--26.
\item[{[B2]}] Beilinson A.A.:{\it  Higher regulators and special
    values of $L$-functions} VINITI, 24  (1984), 181-238. English
  translation: J.Soviet Math., 30 (1985), 2036-2070.
\item[{[BBD]}] Beilinson A.A., Bernstein J., Deligne P.:
  {\it Faisceaux pervers} Asterisque 100, 1982.
\item[{[BD]}] Beilinson A.A., Deligne P.: {\it Motivic polylogarithms and
    Zagier's conjecture}. Preprint
  to appear
\item[{[BMS]}] Beilinson A.A., Macpherson R.D. Schechtman V.V:
  {\it Notes on motivic cohomology}. Duke Math. J., 1987 vol 55
  p. 679-710.
\item[{[BGSV]}] Beilinson A.A., Goncharov A.B., Schechtman
V.V, Varchenko A.N.: {\it Aomoto dilogarithms, mixed Hodge structures and
motivic cohomology of a pair of triangles in the plane }.  The Grothendieck
Feschtrift, Birkhoser, vol , 1990, p 131-172.
\item[{[Bl1]}] Bloch S.: {\it Higher regulators, algebraic $K$-
theory and zeta functions of elliptic curves}, Lect. Notes U.C.
Irvine, 1977.
\item[{[Bl2]}] Bloch S.: {\it Application of the dilogarithm
function in algebraic $K$-theory and algebraic
geometry}. Proc. Int. Symp. Alg. geom., Kyoto (1977), 1-14.
\item[{[Bl3]}] Bloch S.: {\it  Algebraic cycles and higher
    K-theory} Advances in Math.  61, n 3, (1985), 65-79.
\item[{[Bl4]}] Bloch S.: {\it Moving lemma for Higher Chow
    groups} Preprint (1994).
\item[{[Bo]}] B$\ddot o$hm J.: {\it Inhaltsmessung in $\Bbb R^5$ constanter
Kr$\ddot u$mmung },  Arch.Math. 11 (1960), 298-309.
\item[{[Bo1]}] Borel A.: {\it Cohomologie des espaces fibr\'es
principaux}, Ann.\ Math.\ 57 (1953), 115--207.
\item[{[Bo2]}] Borel A.: {\it Cohomologie de $SL_{n}$ et valeurs de
fonctions z\^eta aux points entiers}.  Annali Scuola Normale
Superiore Pisa 4 (1977), 613--636.
\item[{[Ch]}] Chern S.-S.: {\it A simple intrinsic proof of the Gauss-Bonnet
formula for closed Riemannian manifolds}, Ann. of Math., vol 45, 1944,  N 4.

\item[{[Co]}] Coxeter H.S.M. {\it The function of Schlafli and Lobachevsky},
Quart.J.Math. (Oxford) 6 (1935), 13-29
\item[{[Du1]}] Dupont J.: {\it Algebra of polytops and homology of flag
complexes} Osaka J. Math. 19 (1982), 599-641
\item[{[Du2]}] Dupont J.: {\it The dilogarithm as a characteristic
class for a flat bundles}, J.\ Pure and Appl.\ Algebra,
44 (1987), 137--164.
\item[{[DPS]}] Dupont J., Parry W., Sah H.: {\it Homology of classical groups
made discrete.II.$  H_{2}, H_{3}$, and relations with scissors congruences }
J.\ Pure and Appl.\ Algebra, 113,  1988, p.215-260.
\item[{[DS1]}] Dupont J., Sah H.: {\it Scissors congruences II},
J.\ Pure and Appl.\ Algebra, v. 25, 1982, p. 159--195.
\item[{[DS2]}] Dupont J., Sah H.: {\it Homology of Euclidean
    groups of motion made discrete and Euclidean scissor
    congruences},  Acta Math. 164 (1990) 1-24.
\item[{[EP]}] Epstein D.B.A., Penner .R. {\it Euclidean
    decomposition of non-compact hyperbolic manifolds}
  J.Diff. Geom. 27 (1988), 67-80.
\item[{[FV]}] Feigin B, Tsygan B.: {\it Additive K-theory}
  Lect. Notes in Math. 1289.
\item[{[FV]}] Friedlander E., Voevodsky V.: {\it Bivariant cycle
    homology} Preprint 1994.
\item[{[G1]}] Goncharov A.B.: {\it Polylogarithms and motivic Galois
groups} Proseedings of AMS
Reasearch Summer Conference ``Motives'',
Symposium in Pure Mathematics, vol 55, part 2, 43-97.
\item[{[G2]}] Goncharov A.B.: {\it Geometry of configurations, polylogarithms
and motivic cohomology} Advances in Mathematics, vol. 144, N2, 1995. 197-318.
\item[{[G3]}] Goncharov A.B.: {\it Chow polylogarithm and regulators}
Mathematical Research Letters, 1995, vol 2, N 1 95-112.
\item[{[G4]}] Goncharov A.B.: {\it Hyperlogarithms, multiple zeta
    numbers and mixed tate motives}. Preprint MSRI 1993.
\item[{[Goo]}] Goodwillie {\it} Annals of Math. 24 (1986) 344-399.
\item[{[Gu]}] Guichardet A.: {\it Cohomologie des groupes
topologiques et des alg\`ebres de Lie}, Nathan (1980).
\item[{[K1]}] Kellerhals R.: {\it On the volume of hyperbolic 5-orthoschemes
and the trilogarithm }, Preprint MPI/91-68, (Bonn).
\item[{[K2]}] Kellerhals R.: {\it The Dilogarithm and the volume of Hyperbolic
Polytops }, Chapter 14 in [Le].
\item[{[Le]}] {\it Structural Properties of Polylogarithms}, edited by
L.Lewin, AMS series of Mathematical surveys and monographs, vol 37, 1991.
\item[{[Lei]}] Gerhardt C.I. (ed).{\it  G.W.Leibniz Mathematisch Schriften}
III/1, 336-339 Georg Olms Verlag, Hildesheim. New York 1971.
\item[{[Me]}] Meyerhoff R.:  {\it Density of the Chern-Simons invariants for
hyperbolic 3-manifolds}, London Mathematical Sosiety Lecture Notes, vol 112,
1986, 210--226.
\item[{[M1]}] Milnor J.: {\it Computation of volume}, Chapter 7 in [Th].
\item[{[M2]}] Milnor J.: {\it Algebraic $K$-theory and quadratic
forms}, Invent. Math.\ 9 (1970), 318--340.
\item[{[MM]}] Milnor J., Moore J.: {\it On the structure of Hopf
algebras}, Ann. of Math. (2) 81 (1965) 211--264.
\item[{[Mu]}] M$\ddot u$ller P.: {\it  $\ddot U$ber Simplexinhalte in
nichteuklidischen $R\ddot a$umen}, Dissertation Universit$\ddot
a$t Bonn, 1954
\item[{[N]}] Neumann W.: {\it Combinatorics of triangulations and
    the Chern-Simons invariant for hyperbolic manifolds} in
  Topology 90, Proc. of Res. Semesterin Low Dimensional topology
  at Ohio State (Walter de Gruyter Verlag, Berlin - New York
  1992, 243-272
\item[{[N-Z]}] Neumann W., Zagier D.: {\it Volumes of hyperbolic three-
manifolds}, Topology 24 (1985), 307-331.
\item[{[N]}] Neumann W. Yang J.: {\it Invariants from
    triangulations
of hyperbolic 3-manifolds} Electronic Res. Announcements,
    (1995) vol. 1 Issue 2, 72-79
\item[{[P]}] Positselsky L.: {\it Some remarks on triangulated
    categories}. Preprint 1995.
\item[{[Ra]}] Raghunatan M.S. :{\it Discrete subgroups of Lie groups}, Springer
Verlag, 1972.
\item[{[R]}] Ramakrishnan D.: {\it Regulators, algebraic cycles and values of
$L$-functions}, Contemporary Math., vol 83, 183--310.
\item[{[Sah1]}] Sah H.: {\it Hilbert's third problem: scissor
    congruences }Research notes in mathematics 33), Pitman
  Publishing Ltd., San-Francisco-London-Melbourne, 1979
\item[{[Sah2]}] Sah H.: {\it Scissor congruences, I: The Gauss-Bonnet map},
Math. Scand. 49 (1981), 181-210.
\item[{[Sah3]}] Sah H.: {\it Homology of classical groups made
    discrete. III} J.pure appl. algebra 56 (1989) 269-312.
\item[{[SW]}] Sah H., Wagoner J.B: {\it Second homology of Lie groups made
discrete.} Comm. Algebra 5(1977), 611-642.
\item[{[Sch]}]Schl$\ddot a$fli: {\it On the multiple integral
    $\int^{n}dxdy...dz$, whose limits are $p_{1}
= a_{1}x + b_{1}y + ... + h_{1}z > 0, p_{2} > 0, ... , p_{n} > 0$,and $ x^{2} +
y^{2} + ... + z^{2} < 1$}, in Gesammelte Mathematische Abhandlungen II, pp.
219-270, Verlag Birkhauser, Basel, 1953.
\item[{[S1]}] Suslin A.A.: {\it Algebraic $K$-theory of fields}.
Proceedings of the International Congress of Mathematicians,
Berkeley, California, USA, 1986, 222--243.
\item[{[S2]}] Suslin A.A.: {\it $K_{3}$ of a field and Bloch's
group}, Proceedings of the Steklov Institute of Mathematics 1991, Issue 4.
\item[{[Th]}] Thurston W.:{\it Geometry and topology of
    3-manifolds}, Princeton Univercity, 1978.
\item[{[V]}] Voevodsky V.: {\it Triangulated category of motives over a field}
Preprint 1994.
\item[{[Y]}] Yoshida : {\it $\eta$-invariants of hyperbolic 3-manifolds},
Inventiones Math., 1986.
\item[{[Z1]}] Zagier D.: {\it Polylogarithms, Dedekind zeta
functions and the algebraic $K$-theory of fields}. Arithmetic
Algebraic Geometry (G.v.d.Geer, F.Oort, J.Steenbrink, eds.),
Prog. Math., Vol 89, Birkhauser, Boston, 1991, pp. 391--430.
Proceedings of the Texel conference on Arithmetical Algebraic Geometry, 1990.
\item[{[Z2]}] Zagier D.: {\it Hyperbolic manifolds and special
values of Dedekind zeta functions}, Inventiones Math. 83
(1986), 285--301.
\item[{[Z3]}] Zagier D.: {\it The remarkable Dilogarithm}
  J.Math. and Phys. Soc. 22 (1988), 131-145.
\item[{[W]}] Wang H.C.: {\it Topics in totally discontinuous groups}, in
Symmetric spaces, Boothby-Weiss, Editors, New-York, 1972.
\end{itemize}
\end{document}